\pdfoutput=1
\documentclass{amsart}
\usepackage{graphicx}
\usepackage{amssymb,amsmath,latexsym} 
\newcommand{\labell}[1] {\label{#1}}

\newcommand{\vareps}{{\epsilon}}
\newcommand{\Aff}{{\rm Aff}}

\renewcommand{\Tilde}{\widetilde}
\renewcommand{\Hat}{\widehat}
\newcommand{\QED}{\hfill$\Box$}
\newcommand{\less}{{\smallsetminus}}
\newcommand{\NI}{{\noindent}}
\newcommand{\SSS}{{\smallskip}}
\newcommand{\MS}{{\medskip}}

\newcommand{\la}{{\lambda}}
\newcommand{\De}{{\Delta}}

\newcommand{\ga}{{\gamma}}

\newcommand{\Vertt} {{\rm Vert}}
\newcommand\ft {{\mathfrak t}}
\newcommand\p {{\partial}}

\newcommand{\Ee}{{\mathcal E}}

\newcommand{\ov}{\overline}

\newcommand{\Aa}{{\mathcal A}}
\newcommand{\Ss}{{\mathcal S}}
\newcommand{\Cc}{{\mathcal C}}

\newcommand{\La}{{\Lambda}}
\newcommand{\CP}{{\mathbb CP}}
\newcommand{\C}{{\mathbb C}}
\newcommand{\Q}{{\mathbb Q}}
\newcommand{\R}{{\mathbb R}}
\newcommand{\Z}{{\mathbb Z}}

\newcommand{\Om}{{\Omega}}

\newcommand{\om}{{\omega}}
\newcommand{\ka}{{\kappa}}
\newcommand{\al}{{\alpha}}
\newcommand{\si}{{\sigma}}
\newcommand{\io}{{\iota}}
\newcommand{\be}{{\beta}}
\newcommand{\Symp}{{\rm Symp}}

\newcommand{\Ham}{{\rm Ham}}

\newcommand{\Vol}{{\rm Vol}}

\newcommand{\Isom}{{\rm Isom}}

\newcommand{\eps}{{\varepsilon}}

\newtheorem{theorem}{Theorem}[section]
\newtheorem{thm}[theorem]{Theorem}

\newtheorem{cor}[theorem]{Corollary}

\newtheorem{lemma}[theorem]{Lemma}

\newtheorem{prop}[theorem]{Proposition}

\newtheorem{defn}[theorem]{Definition}
\newtheorem{rmk}[theorem]{Remark}
\newtheorem{example}[theorem]{Example}
\newtheorem{quest}[theorem]{Question}
\numberwithin{figure}{section}
\numberwithin{equation}{section}

\title{The topology of toric symplectic manifolds}
\author{Dusa McDuff}
\thanks{Partially supported by NSF grant DMS 0905191}
\address{Mathematics Department, Barnard College, Columbia University
NY, USA}
\email{dmcduff@barnard.edu}
\urladdr{http://www.barnard.edu/mcduff/index.htm}
\date{April 15, 2010, revised September 29 2010}

\begin{document}

\begin{abstract}  
This is a collection of results on the topology of toric symplectic manifolds. Using an idea of Borisov, we  show that a closed symplectic manifold supports at most a  finite number of toric structures.  Further,  the product of two projective spaces of complex dimension at least two (and with a standard product symplectic form) has a unique toric structure.  We then discuss various constructions, using wedging
to build a monotone toric symplectic manifold whose 
center is not the unique point displaceable by probes, and bundles and blow ups  to form manifolds with more than one toric structure.  The bundle construction
uses the McDuff--Tolman concept of  mass linear function.
 Using Timorin's description of the cohomology algebra via the volume function 
 we develop a cohomological criterion for a function to be mass linear, and explain its relation  to Shelukhin's higher codimension  barycenters.
  \end{abstract}
\maketitle

\tableofcontents

%%%%%%%%%%%%%%%%%%%%%%%%%%%%%%%%%%%%%%%%%%%%%%%%%%%%%%%%%%%
\section{Introduction}
%%%%%%%%%%%%%%%%%%%%%%%%%%%%%%%%%%%%%%%%%%%%%%%%%%%%%%%%%%%

The paper \cite{MS} by Masuda and Suh raises many questions 
about the topology of toric manifolds.  One of the most  interesting  can be loosely stated as:

\begin{quest}\labell{qu:1} To what extent does the cohomology 
ring $H^*(M)$  
  determine the toric manifold $M$ or, failing that, the combinatorics of its 
moment polytope?
\end{quest}

Such questions are known under the rubric of {\it cohomological rigidity}; cf. 
Choi--Panov--Suh \cite{CPS}.
One can interpret them in various contexts, including that of 
complex manifolds or
quasitoric (torus) manifolds. In this paper we work exclusively with 
closed symplectic manifolds,  and refine the above question to ask about the 
symplectomorphism type of $(M,\om)$.  Thus our classification is finer than one that considers only the homeomorphism type of $M$ or the combinatorics of the moment polytope, but coarser than one that considers $M$ as a smooth complex variety with
a  given symplectic form.

Recall that a closed symplectic  $2n$-dimensional manifold $(M,\om)$ 
is said to be toric if it supports a Hamiltonian action of an $n$-torus $T$.
This action is generated by a moment map $\Phi:M\to \ft^*$ 
where $\ft^*$ is the dual of the Lie algebra $\ft$ of the
torus $T$.  There is a natural integral lattice $\ft_\Z$ in $\ft$ 
whose elements $H$ exponentiate to circles $\La_H$ in $T$, and hence also
a dual lattice $\ft^*_\Z$ in $\ft^*$.  The image $\Phi(M)$ is well known to be a convex polytope $\De$. It is {\it simple}   ($n$ facets meet at each vertex), {\it rational} (the conormal vectors $\eta_i\in \ft$ to each facet
may be chosen to be primitive and integral),  and {\it smooth} (at each vertex $v$ of $\De$ the conormals to the $n$ facets meeting at $v$ form a basis for
the lattice $\ft_\Z$).  Throughout this paper we only consider such polytopes.   We write them as:
\begin{equation}\labell{eq:De}
\De: = \De(\ka): = \bigl\{\xi\in \ft^*: \langle \eta_i,\xi\rangle \le \ka_i, i=1,\dots,N\bigr\}.
\end{equation}
Thus $\De$ has $N$ facets $F_1,\dots,F_N$ with {\it outward primitive integral conormals} $\eta_i\in \ft_\Z$ and {\it support constants} $\ka=(\ka_1,\dots,\ka_N)\in \R^N$. The  {\it faces} of $\De$ are the intersections $F_I: = \cap_{i\in I} F_i$, where $I\subset \{1,\dots,N\}$. Given a polytope $\De$ we usually denote the corresponding symplectic manifold by $(M_\De,\om_\ka)$.
(See \cite{KKP} for more detailed references on this background material.)

We define $\Cc(\De)$ to be the chamber of $\De=\De(\ka)$, i.e. the open connected  set of all support constants $\ka'$  such that 
$\De(\ka')$ is analogous to $\De(\ka)$; cf. \cite{MT1}. For $\ka,\ka'\in \Cc(\De)$, the symplectic forms $\om_\ka$ and $\om_{\ka'}$ may be joined by the path
$\om_{t\ka + (1-t)\ka'},t\in [0,1],$ and so are deformation equivalent.

Our first result concerns the question of how many different toric actions can be supported by the same symplectic manifold $(M,\om)$.    
Here we identify two toric manifolds 
if there is an equivariant symplectomorphism between them; that is,   
if their moment polytopes may be identified 
by an integral affine transformation.  Karshon--Kessler--Pinsonnault show in \cite{KKP} that  in dimension $2n=4$ a given manifold $(M,\om)$ can support at most a finite number of actions.   The next theorem gives a cohomological version of this result that is valid in all dimensions.
%generalizes  the main theorem in
%Karshon--Kessler--Pinsonnault \cite{KKP} that proves finiteness in dimension $2n=4$.  
Its proof relies on an argument due to Borisov; the original proof applied only when $[\om]$ is integral.  

%The equivariant symplectomorphism type of $(M,\om)$ is well known to be 
%completely determined by the polytope, 

\begin{thm}[Borisov--McDuff]\labell{thm:fin}  Let $R$ be a
commutative ring of finite rank with even grading, 
and write $R_\R: = R\otimes_\Z \R$.
Suppose given
elements  $[\om]\in R_\R$ and 
$c_1, c_{2} \in R$ of degrees $2,2$ and $ 4 $ respectively. 
Then, up to equivariant symplectomorphism, there are at most  
finitely many  toric symplectic manifolds $(M,\om ,T)$ of dimension $2n$
for which there is 
a ring isomorphism $\Psi: H^*(M;\Z)\to R$ that takes the symplectic 
class and the Chern classes $c_i(M), i=1,2,$ 
to the given elements $[\om]\in R_\R, c_i\in R$.
\end{thm}

We prove Theorem \ref{thm:fin} in \S\ref{ss:fin}.  

\begin{rmk}\labell{rmk:fin}\rm (i)
Note that it is crucial to fix the symplectic class $[\om]$ here. Otherwise, 
as is shown by the example of the Hirzebruch surfaces, 
the result is false even for a
 ring as simple as $R=H^*(S^2\times S^2;\Z)$.   In fact,  if $k < \la\le k+1$  the manifold $S^2\times S^2$ with product symplectic form
 $\la pr_1^*\si \oplus pr_2^*(\si)  $  (where $\si$ is an area form on $S^2$)
 supports exactly $k$ different torus actions; cf.  \cite[Example~2.6]{KKP}.
 
\MS

\NI (ii)  One might consider analogous questions for
nontoric symplectic manifolds.  For example, 
one might fix the diffeomorphism type of a closed manifold 
$M$ (rather than its cohomology) 
and fix a cohomology class $a\in H^2(M;\R)$ and ask whether there are 
only finitely many different (i.e. nonsymplectomorphic) symplectic 
structures on $M$ in this class $a$.  The answer here is no:
McDuff \cite{Mcex} constructs an $8$-dimensional manifold that supports infinitely many nondiffeomorphic but cohomologous symplectic forms.    This paper also shows that the manifold $S^2\times S^2\times T^2$ supports infinitely many nonisotopic but cohomologous symplectic forms.  In both cases, the class $[\om]$ is integral
and the forms are deformation equivalent, i.e. they can be joined by a family of (noncohomologous) symplectic forms.    Thus they have the same Chern classes. 
All these examples have nontrivial fundamental group.
Work of Ruan \cite{R} and Fintushel--Stern \cite{FS} shows that in the simply connected case one can find infinitely many nondeformation equivalent symplectic forms on $6$-manifolds of the form $M\times S^2$,
for example when the smooth manifold $M$ is homeomorphic to a $K3$ surface. Although these structures have the same Chern classes, it is not clear whether they can be chosen to be cohomologous.
\MS

\NI (iii)  The extent to which one needs the hypotheses 
on the Chern classes is not clear;
 cf. the discussion in \cite[\S5]{MS}.  By Remark \ref{rmk:finit} they
are unnecessary  if one restricts to
integral $[\om]$.
\MS

\NI (iv)  If one asks the same question in 
the context of $T$-equivariant cohomology, 
then Masuda shows in \cite{Mas2} that the equivariant cohomology
$H^*_T(M;\Z)$, when considered as an algebra over $H^*(BT;\Z)$,
determines the fan, i.e. the family of polytopes 
$\De(\ka), \ka\in \Cc(\De)$,
and hence determines the corresponding toric manifold as a complex variety.
To fix the symplectic manifold, one would also have to specify $\ka$, 
for example by specifying the 
 extension  of the symplectic class to $H^*_T(M;\R)$.
 \MS
 
 \NI (v) Theorem \ref{thm:fin} implies that the number of 
conjugacy classes of $n$-tori in the group $\Ham(M,\om)$ of Hamiltonian symplectomorphisms of $(M,\om)$ is finite, where here we allow conjugation by elements of the full group $\Symp(M,\om)$ of symplectomorphisms of $(M,\om)$. 
Since the orbits of any Hamiltonian action of a torus  are isotropic, each such torus is maximal in $\Ham(M,\om)$.   However, there might be other maximal tori of smaller dimension.  In   \cite{Pin}, Pinsonnault shows that in dimension $2n=4$ there are only finitely many symplectic 
conjugacy classes of such
 maximal tori. 
Again it is important to allow conjugation by elements of $\Symp(M,\om)$;
cf.  \cite[Thm.~1.3]{Pin}.
\end{rmk}

\NI{\bf Manifolds with many toric structures: blow ups and bundles.}
Now  consider the question of which symplectic manifolds 
support more than one toric structure (up to equivariant symplectomrohpism). 
 One easy way to get examples is 
by blowing up points or  other symplectic submanifolds of $(M,\om)$.
(In the combinatorial context the blow up procedure at a point
 is called  vertex cutting; cf.
 \cite[Ex.1.1]{CPS}.)
 
We prove the following result in \S\ref{ss:many}.
%We discuss this  in  \S\ref{ss:many} in the symplectic context. 
Here the weight of a blow up is the symplectic area of the line in the exceptional divisor.   
 
\begin{prop}\labell{prop:torbl}  Suppose that $\De$ is not a product of simplices, and let $(M_\De,\om_\ka)$ be the corresponding symplectic manifold.
Then, for generic choice of $\ka\in \Cc(\De)$,
 there is $\eps_0>0$ such that
any  one point toric blow up $(\Tilde M, \om_{\ka,\eps})$  of $(M_\De,\om_{\ka})$ with weight $\eps<\eps_0$ has at least two  toric structures.
\end{prop}

\begin{rmk}\label{rmk:torbl2}\rm (i) The above result is false when $\De$ is any product of two simplices other than $\De_1\times \De_1$.
This follows because 
 Proposition \ref{prop:uni} below implies that when $\De\ne \De_1\times \De_1$ there is an open nonempty set of $\ka\in \Cc(\De)$ such that 
 $(M_\De,\om_\ka)$ has a unique toric structure, namely that of the product.  Since  
all the vertices of such a product are
equivalent in the sense of Definition \ref{def:equiverr}, its (small)  one point toric blow ups also 
have a unique structure
by Lemma \ref{le:torbl}.\MS

\NI (ii)  To see that one must restrict to generic $\ka$ here, consider the
polygon obtained from  the $2$-simplex
$\De_2$ by blowing up each of  its three vertices
in such a way that all sides of the resulting polygon have equal affine length.
(This polygon corresponds to
 the monotone three point blow up of  
$\C P^2$.)  Then all its vertices are equivalent, so all its one point toric  blow ups are the same.\MS

\NI (iii) We show during the course of the proof of Proposition \ref{prop:torbl} that if the vertices of $\De(\ka)$ are all equivalent for generic $\ka$, then $\De$ is a product of simplices.
\end{rmk}
 
Another natural class of examples is provided by product manifolds 
of the form $M\times S^2$.   The idea is this.  Each $H\in \ft_\Z\less \{0\}$ 
exponentiates to a circle  $\La_H$ in $T$.
Denote by $M_H$ the total space of the
associated Hamiltonian bundle\footnote
{
A smooth bundle $E\to B$ with fiber $F$ is Hamiltonian if its structural group
reduces to the Hamiltonian group $\Ham(F,\si)$ of
some symplectic form $\si$ on $F$. Often, as here,
$\si$ is given.  When $\pi_1(B) = 0$ this is equivalent to saying that 
the fiberwise symplectic form $\si$ extends to a closed form on the 
total space.}
$(M,\om)\to M_H\stackrel{\pi}\to S^2$. 
One can realize $M_H$ as the quotient $S^3\times _{S^1}M$ where 
$S^1$ acts diagonally on $S^3\subset \C ^2$ and via $\La_H$ on $M$.
Consider the $1$-form 
$$
\al: = \tfrac i{4\pi} \Bigl(\sum_{j=1,2} z_j\,d\ov z_j - \ov z_j \,dz_j\Bigr),
$$
on $S^3$.  (The form $\al$ is the standard contact form normalized so that the integral
of $d\al$ over the unit disc $\{(z_1,0): |z_1|<1\}$  is $1$.)   Then the form
$pr^*(\om) +d\bigl((\la-H)\al\bigr) $, where $\la\in \R$ and  
$pr: S^3\times M\to M$ is the projection, descends to the quotient
 $M_H$ and defines a symplectic form  $\Om_\la$ there provided that the function $\la-H$ is positive on $M$.
%Alternatively
%(CHECK SIGNS!)
%$$
%M_H = (\bigl(D_+ \times M\bigr)\cup\bigl(D_-\times M\bigr)/\!\sim,\quad
%(e^{2\pi it},\phi_tx)_+ \equiv (e^{2\pi it},x)_-,
%$$
%where $\la_H = \{\phi_t\}_{t\in [0,1]}$, and $D_\pm$ are unit $2$-discs.
% The fiberwise 
%symplectic form extends to a closed form $\Om$ on $M_H$ that 
%we normalize by requiring that
%$\int _{M_H}\Om^{n+1} = 0$.  It is easy to see that for suitably large 
%constant $\la>0$ the form $\Om_\la: = \la \pi^*(\al)+\Om$ 
%on $M_H$ is symplectic, where 
%$\al$ is the $S^1$ invariant symplectic form on $S^2$.
Moreover, because the action of $\La_H$ on $M$ commutes with $T$
the manifold
$(M_H,\Om_\la)$ supports an action of $T_H: = T^{n+1}$.  
Thus $(M_H,\Om_\la)$ is toric.  Moreover the bundle 
\begin{equation}\labell{eq:torbun}
M\stackrel{\io}\hookrightarrow M_H
\stackrel{\pi} \to S^2
\end{equation}
 is toric in the sense that there is a 
group homomorphism $\rho:T_H\to S^1$ such that the projection 
$\pi:M_H\to S^2$
intertwines the action of $T_H$ on $M_H$ with the action of
 $\rho(T_H)=S^1$ on $S^2$; i.e.
$$
\pi(t\cdot x) = \rho(t)\cdot \pi(x),\quad x\in M_H,\;t\in T_H.
$$

Suppose now that $\La_H$, when considered 
as a loop in the Hamiltonian group $\Ham(M,\om)$, is contractible.
Then the bundle $(M,\om) \to M_H\to S^2$ is trivial as a Hamiltonian bundle.  
This readily implies\footnote
{
for example by adapting the proof of Proposition 9.7.2 (ii) on p 341 of \cite{JHOL}.}
  that $(M_H,\Om_\la)$ is symplectomorphic to the product
$\bigl(S^2\times M,\si \oplus \om\bigr)$ for suitable area form $\si$ on $S^2$. 
But we will see in Remark \ref{rmk:bun} (i)
 that the moment polytope $\De_H$ of $(M_H,\Om_H,T_H)$ is not affine equivalent to a product when $H\ne 0$.   
This proves the following result.
 
\begin{lemma} \labell{le:contr} If $(M,\om,T)$ is such that  
the loop $\La_H$ contracts in $\Ham(M,\om)$ for some  nonzero $H\in \ft$
then there is $\la_0>0$ such that  for all $\la\ge \la_0$ the product 
$\bigl(S^2\times M,\la\si \oplus \om\bigr)$ supports
more than one toric structure.
\end{lemma}

\begin{rmk}\labell{rmk:contr}\rm (i)  The moment polytope $\De_H$ of $M_H$ is
always
combinatorially equivalent to a product.  
Hence the examples in Lemma \ref{le:contr} are not distinguished 
in papers such as \cite{CPS}.  \MS

\NI (ii)  One could, of course, also consider the (toric) 
bundle $S^{2k+1}\times _{S^1}M\to \C P^k$ corresponding to
 the loop $\La_H$ for $k>1$.  However, even if 
 $\La_H$ contracts in $\Ham(M,\om)$,
  this bundle is {\it never} trivial as 
a Hamiltonian bundle when $H\ne 0$; cf. Remark \ref{rmk:trivk}.
\end{rmk}

The next question is: when do such loops exist?  The paper
 McDuff--Tolman \cite{MT1} analyses this question in great detail.  
The easiest case is when the loop $\La_H$ (or one of its finite multiples 
$\La_{mH}$) contracts in the maximal 
compact subgroup $\Isom_0(M)$ of $\Ham(M,\om)$, consisting of 
symplectomorphisms that preserve the natural K\"ahler metric
on $(M,\om)$.\footnote
{
This group is described in slightly different language in Masuda \cite{Mas}.}  
Such elements $H\in \ft_\Z$ were called {\bf inessential} in \cite{MT1},
and exist when the moment polytope $\De$ of $M$ 
satisfies some very natural geometric conditions.
In particular, by \cite[Prop.~3.17]{MT1} if they exist the
polytope $\De$ must either be a bundle over a simplex or an expansion 
(wedge).   Correspondingly $M$ is either the total space of a 
toric bundle over $\C P^k$ or is the total space of smooth Lefschetz pencil
with axis of (real) codimension $4$. (The last statement is explained in
\cite[Rmk.~5.4]{MT1}.)  Generalized Bott towers, which are iterated bundles formed from projective spaces, are well known examples.
Since any wedge and any bundle over $\C P^k$  has a nontrivial inessential function $H$,  many product toric manifolds $M\times S^2$ 
have more than one toric structure. For further discussion of 
this issue, see Theorems \ref{thm:coh} and \ref{thm:FML2} below.
\MS

\NI{\bf Manifolds with unique toric structures.}
Next, one might wonder which symplectic manifolds have just one toric structure.  We prove the following result in \S\ref{ss:uni} by a 
cohomological argument.   We denote by $\om_n$ the usual  symplectic form on $\C P^n $   that integrates over a line to $1$.  Thus $(\C P^n,\om_n)$ is a toric manifold with moment polytope 
equal to the standard unit simplex 
$$
\De_n=\bigl\{x_1\ge 0,\dots,x_n\ge 0, \sum x_i \le 1\bigr\}\subset \R^n.
$$ 

\begin{prop}\labell{prop:uni}
Let $(M,\om) = \bigl(\C P^{k}\times \C P^{m},\om_{k} + \la\om_{m}\bigr)$, where  $\la > 0$.  If $k\ge m\ge 2$  then $(M,\om)$ has a unique toric structure.  If $k>m=1$, this remains true provided that
 $\la \le1$, while if $k=m=1$ we require $\la = 1$.
\end{prop}

\begin{rmk}\rm (i)
At first glance, this result  is somewhat surprising, since one might well 
imagine that there are analogs of  Hirzebruch structures
 on products such as $\C P^3\times \C P^2$.   As pointed out in Remark \ref{rmk:trivk}, the explanation for this lies in the characteristic classes constructed in \cite{KM}.\MS
 
 \NI (ii)
If $k\ge m=1$ and $\la>1$, then there are nontrivial
 toric $\CP ^k$ bundles over $\CP^1$   that are  symplectomorphic to products
for large $\la$, as one can see by arguments similar to those 
that prove Lemma \ref{le:contr}.  However, even in the case $k=m=1$,
when we get the  Hirzebruch surfaces, the proof that these manifolds are symplectomorphic to products for all relevant $\la$ is nontrivial; see \cite{Mcex} or \cite[Prop.~9.7.2]{JHOL}.   Nevertheless, this proof
 should generalize to show that uniqueness fails whenever $\la$ does not 
 satisfy the conditions in Proposition \ref{prop:uni}.
 \MS
 
 \NI (iii) Proposition \ref{prop:uni} extends work by Choi, Masuda and Suh, who show in \cite{CMS} that if $M$ is a toric $\C P^k$-bundle over 
 $\C P^m$ then it is diffeomorphic to the product of its base and fiber exactly if its integral cohomology ring is isomorphic to that of the product.
\end{rmk}

\NI {\bf Monotone polytopes.}\,\,
Another natural class of manifolds that might have unique toric structures  is that of
 monotone  manifolds.  Recall that a
 %It is convenient to normalize the scale of the symplectic form.
 %Thus %We will slightly change the usual definition as follows: 
 symplectic manifold $(M,\om)$ is  said to be
{\bf monotone} if  there is $\la>0$ such that $[\om] =  \la c_1(M)$.  In this paper, we shall always normalize $\om$ so that $\la = 1$.  Thus, in the toric case, the 
moment polytope is scaled so that the affine length of each edge $\vareps$ 
is precisely $\int_{\Phi^{-1}(\vareps)} \om$.  

The moment image of a monotone toric manifold is called a monotone polytope.\footnote
{
These are also known as smooth reflexive   polytopes.  Note that much of the literature about them is written in terms of their dual polytopes $P\subset\ft$
(which are simplicial) rather than the moment polytopes considered here.} 
 Since rather little seems to be known in general about their structure,  we 
 begin our discussion by  describing some 
 elementary constructions.  
% For example, in \S\ref{ss:cut} we give a \lq\lq 
%new" $4$-dimensional shape (due to Paffenholz) obtained by cutting (a 
%generalized blow up). 
%  It is new in the sense that it is not a wedge or a  bundle formed from lower dimensional polytopes or formed from such by blow up. 
 
 The most interesting of these   is that of wedge (called expansion in \cite{MT1}).  It was used in Haase--Melnikov \cite{HM} to show that {\it every} smooth integral polytope is the face of some monotone polytope.
 We  adapt it here to answer some questions raised in \cite{MT1}.
  Let us say that a facet  $F$ of a polytope is {\bf pervasive} if it meets all other facets and is {\bf powerful} if there is a edge between $F$ and every vertex of $\De$ not on $F$.   We showed in \cite[Thm.~A.6]{MT1} that 
 in  dimension $\le 4$ the only  polytopes with all  facets  powerful are 
 combinatorially equivalent to products of simplices. This is not true in higher dimensions, even if one restricts to the monotone case.    For every face of a product of simplices is also a product of simplices, while, by   Lemma \ref{le:pp},  a monotone polytope
 has faces of arbitrary shape.
 
 The next result is an immediate consequence of 
  Lemma \ref{le:pp}.
  
  \begin{prop} \labell{prop:pp} Let $\De'$ be any smooth polytope with integral vertices.  Then some multiple $k\De', k\in \Z,$ is integrally affine equivalent to a face in a monotone polytope $\De$  all of whose facets are pervasive and powerful.
  \end{prop}

  Further, in Lemma \ref{le:sE} we  use the wedge construction
   to describe an example found by Paffenholz of a monotone polytope 
  that fails the star-Ewald condition of \cite{Mcpr}.  
  As we explain in \S\ref{ss:wedge}, this is related to the work of Fukaya--Oh--Ohta--Ono \cite{FOOO} on the Floer homology of toric fibers.  

Despite the existence of this rather versatile construction, I do not know the answer to the following question.

\begin{quest}\labell{qu:3} Is there a monotone toric 
manifold $(M,\om)$ with more than one toric structure?
\end{quest}

It is not clear whether one can obtain such an 
example by blowing up a point (vertex cutting).   However, the
 next result shows that one cannot get examples by the bundle construction 
used in Lemma \ref{le:contr} above. 

We shall say that two bundles 
$M\to M_{H_i}\to S^2, i=1,2,$ are bundle isomorphic 
if there is a commutative diagram
$$
\begin{array}{ccccc} M&\to &M_{H_1}&\to& S^2\\
id\downarrow&&\phi\downarrow&&  id\downarrow\\
M&\to &M_{H_2}&\to& S^2\end{array}
$$
where $\phi$ is a diffeomorphism.  Thus we assume that $\phi$ is the identity 
map on the distinguished fiber.  However, it need not preserve the
symplectic forms on the total spaces.   

\begin{defn}\labell{def:equiv}  We say that two facets $F_i, F_j$ of $\De$ are {\bf equivalent}, and write 
%equivalence relation on the facets of a polytope $\De$ under which 
$F_i\sim F_j$, if there is a vector $\xi\in \ft^*$ that is parallel to 
all other facets of $\De$.
\end{defn}

  It is shown in 
\cite[Lemma~3.4]{MT1} that  $F_i\sim F_j$ precisely if there is a robust\footnote
{This means that the reflection persists as one perturbs $\ka$ a little.}
 affine reflection of $\ft^*$ that takes $\De$
 to itself and interchanges the facets $F_i, F_j$,  fixing all others.
Because it is robust, this affine reflection lifts to a symplectomorphism 
of $(M_\De,\om_\De)$
that
lies in the maximal compact subgroup $\Isom_0(M_\De)$
of $\Ham(M_\De,\om)$; in particular it is isotopic to the identity.  
It also follows from the Stanley--Reisner presentation of $H^*(M)$ (cf. equation \eqref{eq:SR})
 that $F_i\sim F_j$ exactly if the hypersurfaces 
$\Phi^{-1}(F_i)$ and $\Phi^{-1}(F_j)$ represent the same 
element in $H_{2n-2}(M)$.

We prove the following result in \S\ref{ss:bun}.
 
\begin{prop}\labell{prop:monprod}  Suppose that
$(M_H, \om_H,T_H)$ is a monotone toric manifold with 
moment polytope $\De_H$
that is the 
total space of a toric bundle with fiber $(M,\om,T)$ and base $\C P^1$. 
%as in equation \eqref{eq:torbun}. 
Then the following hold.
\begin{itemize}\item[(i)] Either there is a facet $F_j$  of the moment 
polytope $\De$ of $M$ such that $H = \eta_j$, or 
$H=0$ and $\De_H$ is affine equivalent to
the  product $\De_1\times \De$.
\item[(ii)] If $H=\eta_j$ the loop $\La_{H}$ does not contract
in $\pi_1\bigl(\Ham(M_\De,\om_\De)\bigr)$, and $(M_H,\om_H)$ is not symplectomorphic to a product $(M\times S^2,\om \oplus \si)$.
\item[(iii)] Two of the bundles in (ii)  are bundle isomorphic only if
they are generated by elements $H_j = \eta_j, j=1,2,$ that correspond to 
equivalent facets of $\De$. In this case, the loops 
$\La_{H_j}$ are conjugate in $\Ham(M_\De,\om_\De)$.
\end{itemize} 
%Suppose that the monotone toric manifold $(M,\om_M,T)$ is a toric 
%bundle over $S^2$ with fiber the toric manifold $(M_{\Tilde\De}, 
%\om,\Tilde T)$.   Then for each possible fiber
%$(M_{\Tilde\De}, \om,\Tilde T)$ there is at most one such structure, up to fiberwise symplectomorphism.
\end{prop} 

\MS
\NI
{\bf Mass Linearity.}
Our final set of results again concerns the question of which toric manifolds
  $(M,\om)$ have nontrivial loops $\La_H$ that contract in $\Ham(M,\om)$.
Above we discussed inessential $H$.\footnote
{
By slight abuse of language, we often call $H$  a function, thinking of it as 
a function on the moment polytope $\De$. Note also that the moment map for the 
circle action $\La_H$  is the composite
$x\mapsto \langle H,\Phi(x)\rangle$, of $\Phi$ with the projection $\ft^*\to \R$ given by  inner product with $H$.}
  The papers \cite{MT1,MT2}
discuss a more 
interesting class of functions 
$H$ called {\bf mass linear functions}.  These are functions on $\De$ 
whose value $H(B_n)$ at the barycenter $B_n(\ka)$ of the moment polytope
$\De=\De(\ka)$ is a linear function of
the support numbers $\ka = (\ka_1,\dots,\ka_N)$ of its facets.
By \cite[Prop.~1.17]{MT1} every inessential function is mass linear. 
However, even when $n=3$ there are pairs $(\De,H)$ where $H$ is mass linear but is not inessential; in this case we say that 
$H$ is {\bf essential}. 
By \cite[Thm.~1.4]{MT1}, in $3$ dimensions there is 
precisely one such family  $(\De,H)$ 
that we describe in Lemma \ref{le:prod} below.   In these examples, 
the underlying polytope 
$\De$ is a $\De_2$-bundle over $\De_1$, where $\De_k$ denotes the
standard  $k$-simplex.\footnote
{
Unless explicit mention is made to the contrary, we allow the standard simplex to have any size, i.e. we do not fix $\ka$.}

We showed in  \cite[Prop.~1.22]{MT1} that
 if a loop $\La_H$ contracts in $\Ham(M,\om)$ then
 $H$ is mass linear.  There the  argument  was based on
 Weinstein's action homomorphism of $\pi_1(\Ham(M,\om)$; in \S\ref{ss:shel} below 
 we explain an alternative argument
 due to Shelukhin that uses some other homomorphisms.
 Conversely, one can ask if the  mass linearity of $H$ implies that
 the loop $\La_{mH}$ contracts in $\Ham(M,\om)$ for some $m$.
 (Proof that this is true in some nontrivial cases
is the subject of ongoing research.)
 Our next result establishes a cohomological version 
 of this statement.

We prove the following result in
\S\ref{s:FML}, using Timorin's very interesting   
description of the real cohomology algebra  of $(M,\om_\ka)$ 
in terms of the 
function $V(\ka)$ that gives the volume of the moment polytope 
in terms of the support numbers $\ka$.
%to prove the following result.

\begin{thm}\labell{thm:coh} 
Let $(M,\om,T)$ be a toric manifold with moment polytope $\De$, and let 
$H\in \ft\less \{0\}$.  Let $M\to M_H\to S^2$ be the corresponding bundle.
Then the element $H\in \ft_\Z$ is mass linear 
if and only if
there is an algebra isomorphism 
$$
\Psi:
H^*(S^2;\Q)\otimes H^*(M;\Q) \stackrel{\cong}\to
H^*(M_H;\Q) 
$$
that is compatible with the fibration structure in the sense that it fits into a commutative diagram
$$
\begin{array}{ccccc} H^*(M)&\leftarrow &H^*(S^2)\otimes H^*(M)&\leftarrow& H^*(S^2)\\
id\downarrow&&\Psi\downarrow&&  id\downarrow\\
H^*(M)&\leftarrow & H^*(M_{H})&\leftarrow& H^*(S^2).\end{array}
$$
\end{thm}

\begin{rmk}\rm If one writes $\Psi$ in terms of a basis for the integral cohomology, then its coefficients give information about the 
 order of the loop $\La_H$ in $\pi_1\bigl(\Ham(M,\om)\bigr)$.  Indeed,
 if this order is $m<\infty$ then these coefficients must lie in $\frac 1m \Z$;
cf. Remark \ref{rmk:coeff}.
%
%\NI (ii)
%If we could find  an essential mass linear $H$ for 
%which some multiple of $\La_H$ contracts in 
%$\Ham(M,\om)$,  we would by
%Lemma \ref{le:contr}  have a nontrivial 
%example where the (integral)  cohomology ring does 
%determine the toric structure.
\end{rmk}

Theorem \ref{thm:cohom} below sharpens Theorem \ref{thm:coh}, using  Shelukhin's
concept of full mass linearity. 
He considers all the barycenters $B_k$, 
 $k=0,\dots,n,$ of $\De$, defining $B_k$ to be
the barycenter of the union of the
 $k$-dimensional faces of $\De$. For example, $B_0$ is the 
average of the vertices of $\De$.
He showed that the numbers $H(B_k)-H(B_n)$ are the values 
of some natural characteristic classes on toric loops $\La_H$,
hence 
proving the following result.

\begin{prop}[Shelukhin \cite{Shel}]
The loop $\La_H$ contracts in $\Ham(M_\De,\om_\ka)$ only if 
$H(B_n) = H(B_k)$ for all $k=0,\dots,n-1$.
\end{prop}

We will say that $H$ is {\bf fully mass linear} if $H(B_k) = H(B_n)$ 
for  $0\le k\le n-1$. 
Theorem \ref{thm:cohom} gives a cohomological interpretation of the 
 full mass linearity condition.  In \S\ref{s:FML} we also sharpen 
some of the combinatorial results of \cite{MT1}, obtaining the 
following results.

\begin{thm}\labell{thm:FML2} {\rm (i)} An element $H\in \ft_\Z$ is mass linear
if and only if $H(B_n) = H(B_0)$.  Moreover, in this case, $H(B_{n-1}) = H(B_n).$
\MS

\NI {\rm (ii)}  Every mass linear function on a polytope of dimension at most $3$ 
 is fully mass linear.
\end{thm}

 \NI {\bf Organization of the paper.}
 We begin the proofs by  discussing the structure of monotone manifolds, since this will allow us to introduce some of the main constructions.  
Theorem \ref{thm:fin} is proved in \S\ref{ss:fin}; the argument does not use \S\ref{s:mono}. Mass linearity is discussed in \S\ref{s:FML}.  This section is essentially independent
of the other two. \MS
 
 \NI {\bf Acknowledgements.} This paper owes much to Lev Borisov who
 sharpened the original version of the finiteness theorem so that it applies to all symplectic manifolds, not just to those with integral symplectic form.
   I also warmly  thank Sue Tolman, 
Yael Karshon, Megumi Harada, Tara Holm and Olga Buse  for interesting discussions, Mikiya Masuda for some useful comments,
 and Andreas Paffenholz for using {\it Polymake} to investigate 
the star-Ewald condition.    I am very grateful to 
Seonjeong Park for pointing out some inaccuracies, and 
also to the referee for pointing  out
 many small   mistakes and unclear explanations, which helped  me improve the paper significantly. 
 The NSF-supported 
 Great Lakes Geometry Conference in Madison, April 2010, provided an excellent forum in which  to discuss 
 the results.  This paper was written during my stay at MSRI in Spring 2010, and I am very grateful for its hospitality and support.

%%%%%%%%%%%%%%%%%%%%%%%%%%%%
  \section{Monotone polytopes}\labell{s:mono}
 %%%%%%%%%%%%%%%%%%%%%%%%%%%%

We begin with a general remark about normalizations.  The moment polytope $\De\subset \ft^*\cong \R^n$ of 
a toric manifold $(M,\om,T)$ is determined as a subset of $\R^n$ up to the action of 
the integral affine group $\Aff(n;\Z)$.   Because the conormals at any vertex form a lattice basis, we may therefore always choose coordinates on $\R^n$ so that  the 
conormals at any chosen vertex $v$ are $-e_1,\dots,-e_n$, i.e. the negatives of the standard basis vectors.   Then the polytope lies in a translate of the positive quadrant 
$x_i\ge 0, i=1,\dots,n$.  Sometimes we normalize so that $v=0$, but often (as in the monotone case considered below) we set $v=(-1,\dots,-1)$ so that the center point  of $\De$ is at $\{0\}$.

Recall that the symplectic manifold $(M,\om)$ is monotone if $[\om] = \la c_1(M)$ for some $\la>0$.  Throughout we will normalize monotone manifolds  so that $\la =1$.  
  There are several possible ways of characterizing the moment image of a monotone toric manifold.  
%  definitions of a monotone (moment) polytope.
The following  well-known lemma is proved in \cite[Lemma~3.3]{Mcpr}.

\begin{lemma}\labell{le:mono}
 A  simple smooth polytope $\De$  is monotone if and only if it satisfies the following conditions:
 
\begin{itemize}\item[(i)]
$\De$ is an integral  (or lattice) polytope in $\R^n$ 
with a unique interior integral point  $u_0$,

%\NI $\bullet$ it  satisfies the Delzant smoothness condition: at each vertex $v$ the $n$ primitive integral vectors  $e_i(v)$ pointing along the edges from $v$ form a basis for the integral lattice $\Z^n$,\SSS

\item[(ii)]  $\De$ satisfies the  {\bf vertex-Fano} condition:
for each vertex 
$v_j$ we have
$$
v_j+\sum_ie_{ij} = u_0,
$$ 
where  $e_{ij}, 1\le i\le n,$ 
are the primitive integral
 vectors from $v_j$ pointing along the edges of $\De$.
 \end{itemize}
\end{lemma}

\begin{rmk}\rm (i)  If the conditions in Lemma \ref{le:mono} are satisfied,  then
the affine distance\footnote
 {
 See \cite[\S2]{Mcpr} for a general explanation of how to measure affine distance.} 
  $\ell_j(u_0): = \ka_j - \langle \eta_j,u_0\rangle$ from $u_0$ to the facet $F_j$ equals $1$  for all $j$.   Hence if  we translate $\De$ so that $u_0=\{0\}$ the structure 
 constants $\ka_i$ in the formula (\ref{eq:De}) are all equal to $1$. 
  Conversely, any integral polytope with $\ka_i=1$ for all $i$
   satisfies conditions (i) and (ii) in Lemma \ref{le:mono} with $u_0 = \{0\}$ and so is monotone.  \MS
% for every vertex  $v$ one can choose coordinates for which $u_0=(0,\dots,0)$, $v=(-1,\dots,-1)$  and the facets through $v$ are $\{x_i=-1\}$, $i=1,\dots,n$.  In particular the affine distance\footnote
% {
% See \cite[\S2]{Mcpr} for an explanation of how to measure affine distance.} 
%  $\ell_j(u_0): = \ka_j - \langle \eta_j,u_0\rangle$ from $u_0$ to the facet $F_j$ equals $1$  for all $j$.  Thus if we translate $\De$ so that $u_0=\{0\}$ the structure 
% constants $\ka_i$ in the formula (\ref{eq:De}) are all equal to $1$. 
%  Conversely, any integral polytope with the properties described in the last sentence satisfies conditions (i) and (ii) in Lemma \ref{le:mono}.  \MS

 \NI (ii) Another closely related notion is that of   Fano polytope.
Usually one defines this  in terms of the dual $P\subset \ft$ to the moment polytope (namely the fan), and calls $P$  Fano  if one can choose
support constants $\ka'$ for the moment polytope $\De$ that make it monotone.   However, the constants $\ka'$ are not specified.
Correspondingly, a {\bf Fano toric symplectic  manifold} $(M,\om_\ka,T)$ is one that may not be monotone but where there is $\ka'\in \Cc_\De$ such that
$(M,\om_{\ka'},T)$ is monotone.
\end{rmk}

 %%%%%%%%%%%%%%%%%%%%%%%%%%%%
  \subsection{The wedge construction}\labell{ss:wedge}
 %%%%%%%%%%%%%%%%%%%%%%%%%%%%
  
  This is a very useful construction that appeared in \cite{MT1}
because  of our result that any polytope with a nontrivial robust\footnote
{
A  nontrivial  affine transformation of $\De(\ka)$ is called robust if it   persists when one perturbs $\ka$; for a more precise definition see  \cite[Def.~1.11]{MT1}. These symmetries make up the group $\Aff_0(\De)$ of Definition \ref{def:equiverr} below.} 
symmetry
    is either a bundle over a simplex or is an expansion; cf. \cite[Prop.~3.15]{MT1}. Moreover, a polytope has such a symmetry exactly if the identity component of its K\"ahler isometry group (with respect to the natural 
    K\"ahler metric)
      is larger than the torus $T^n$; cf. \cite[Prop~5.5]{MT1}.  We called this construction an {\it expansion}.  However, it is known in the combinatorial literature as a {\it wedge. }

  Here is the definition.  
    
  \begin{defn}\labell{def:expan}  Suppose that $\De\subset \R^n$ is described by the inequalities
 \begin{equation}\labell{eq:poly}
  \langle \eta_i, x\rangle \le \ka_i, \quad x\in \R^n,\;\; i\in \{1,\dots,N\},
\end{equation}
  where
   $\ka_i>0$ so that  $\{0\}$ lies in its interior.  Its {\bf wedge} (or {\bf expansion})
   $\De'$ along the facet $F_k$ lies in $\R^{n+1} = \R^n\times \R$ and is given by the above inequalities for $i\ne k$ (where we identify $\eta_i\in \R^n$ with $(\eta_i,0)\in \R^{n+1}$) together with
  $$
  x_{n+1}\ge -1, \quad \langle \eta_k, x\rangle + x_{n+1}\le \ka_k -1.
  $$
  \end{defn}
  
    Thus  we replace the conormal $\eta_k$ by the two conormals 
    $\eta_k' = (\eta_k,1)$ and $\eta_{N+1}' = (0,\dots,0,-1)$.  
    The original polytope $\De$ is now the facet $F_{N+1}'$ of 
     the wedge $\De'$.  In fact,
    $\De'$  is made from the product $\De\times [-1,\infty)$  by adding a new \lq\lq top" facet $F_k'$  with conormal $\eta_k'=(1,\eta_k)$  that intersects 
  the \lq\lq bottom" facet  $F_{N+1}': = \{x_{n+1}=-1\}$ in the facet $F_k$ of $\De$.  The corresponding toric manifold $M_{\De'}$  is the total space of a smooth Lefschetz pencil with pages $M_\De$ and axis 
  (of complex codimension $2$)  $F_k'\cap F_{N+1}\cong F_k$; cf \cite[Rmk.~5.4]{MT1}.

Note that all the structural constants $\ka_j$ remain the same, 
except for $\ka_k$ which decreases by $1$.  Moreover $\ka_{N+1} = 1$.
Haase and Melnikov point out in \cite[Prop.~2.2]{HM} that
     by repeating this construction until each $\ka_j = 1$ one 
     finds that every integral polytope with an interior integral point (which we can assume to be at $\{0\}$)  is integrally affine equivalent to the face of some monotone polytope.
Here is a slight refinement of their result. Recall that a facet $F$ is called {\it pervasive} if it meets all other facets and {\it powerful} if there is a edge between $F$ and every vertex of $\De$ not on $F$.

\begin{lemma}\labell{le:pp}  Suppose that $\De$ is a smooth integral polytope  with $\{0\}$ in its interior and with all structural constants $\ka_i\ge 2$.
Then $\De$ is a face in a monotone polytope for which all facets are both pervasive and powerful.
\end{lemma}
\begin{proof}
The new facets $F_{N+1}$ (the bottom) and $F_k'$ (the top) of any wedge are pervasive.  Moreover, any pervasive facet of $\De$ remains pervasive in $\De'$.  
Similar remarks apply to the concept of powerful since all vertices in $\De'$ lie either on the top or bottom facet of $\De'$.
The hypothesis that $\ka_i\ge 2$ implies that we must wedge at least once along each facet to get a monotone polytope.
The result follows.
\end{proof}
  
In \cite{MT1} we were interested in polytopes for which all facets are both pervasive and powerful because we were trying to understand mass linear functions  $H$ on polytopes $\De$.  Our basic question was:  is it always true that after subtracting an inessential function $H_0$, the resulting mass linear function has a symmetric facet?\footnote
{
A facet $F_j$ is called {\it symmetric} (resp. {\it asymmetric})  if, when we write
 $H(B_n) = \sum \ga_i\ka_i$, the coefficient $\ga_j$ vanishes (resp. $\ga_j\ne 0$).}
  Equivalently, is there an inessential $H_0$ such that $H-H_0 = \sum \ga_i\ka_i$ where $\ga_i=0$ for some $i$?
  The answer would be yes, if every polytope with all facets powerful and pervasive has at least two equivalent facets; cf. \cite[Lemma~3.19]{MT1}.
Therefore,  it would be relevant to know the answer to the following question.

\begin{quest}  Is there a smooth polytope whose  facets are powerful and pervasive and have the property that no two facets are equivalent?
\end{quest}

Of course, to construct such a polytope one cannot use wedging, since the top and bottom facets of a wedge are always equivalent.

We end this subsection by using wedges to construct an example of a smooth monotone polytope $\De$  that does not satisfy 
the star-Ewald condition of \cite[Definition~3.5]{Mcpr}  at one of its vertices.  
This is a condition on each face $f$ of $\De$ that is designed so that  it fails at $f$
exactly if there is a point in the interior of the cone $C(f,0)$ spanned by $f$ and $\{0\}$ that cannot be displaced by a probe; cf. the proof of \cite[Theorem~1.2]{Mcpr}.
Therefore the corresponding Lagrangian toric fibers $L(u)$ in $M_\De$  may perhaps be nondisplaceable by Hamiltonian isotopies, even though, according to \cite{FOOO}, their Floer homology vanishes.

\begin{figure}[htbp] %  figure placement: here, top, bottom, or page
   \centering
  \includegraphics[width=3in]{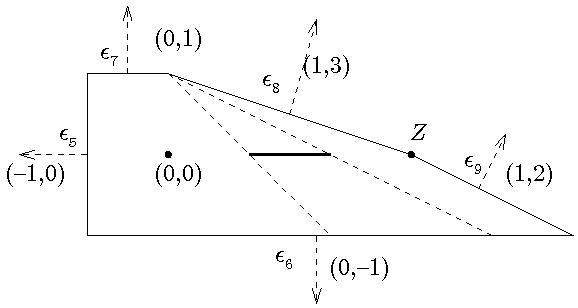} 
   \caption{The polygon $\Tilde\De$; the heavy line segment in the middle
   are the points that are nondisplaceable by probes. The conormals  of
   $\Tilde\De$ are the vectors $\nu_5,\dots,\nu_9$ in equation \eqref{eq:sE}
  with $\ka= (1,1,1,3,3)$.}
   \label{fig:1}
\end{figure}

 This example is due to Paffenholz \cite{Paf}.  
  By using the program {\it Polymake} 
 he shows that all polytopes of dimensions less than $6$ do satisfy the star-Ewald condition.  However, he found three  $6$-dimensional examples  where the condition fails, and many more $7$-dimensional ones. In all but one case the condition failed at a vertex or an edge, but there is one
 $7$-dimensional example (se.7d.02 on his list) where it fails on a nonconvex set consisting of two edges.\footnote
 {If you look at the file, the first edge  together with its data is listed first, 
 and the data on the second edge occurs about half way through.}
 All of his examples are wedges.
  We shall explain the easiest one, which is a repeated wedge of the polygon
  in Figure \ref{fig:1}.   
  
Consider the monotone $6$-dimensional polytope $\De$  with conormals:
 \begin{gather*} \eta_i=-e_i, i=1,\dots, 6, \;\; \eta_7 = e_6,
 \eta_8 = (1,1,0,0,1,3),\;\;\eta_9=(0,0,1,1,1,2),
 \end{gather*}
 where $e_i$ are the standard basis in $\R^6$ and we set $\ka_i=1$ for all $i$.
   Further
 let $\Tilde\De$ be the polygon with conormals
\begin{equation}\labell{eq:sE}
  \nu_5=(-1,0),\;\;\nu_6= (0,-1),\;\;\nu_7=(0,1), \;\;\nu_8=(1,3),\;\;\nu_9 = (1,2),
\end{equation}
  and  with $\ka=(1,1,1,3,3)$ as in Figure \ref{fig:1}.
  Then $\Tilde\De$ can be identified with the facet $F_{0123}$ of $\De$.
Further $\De$ is obtained from $\Tilde\De$  by making twice repeated  expansions in the edges $\vareps_8: = F_ {12348}, \vareps_9: = F_{12349}$ of $\Tilde\De$ (or, more precisely, in the facets corresponding to these edges).
The facets $F_1,F_2$ come from the expansion along $\vareps_8$ and the facets $F_3,F_4$ from the expansion along  $\vareps_9$. 

Given an integral polytope $\De$ with $\{0\}$ in its interior, 
consider the set  
$$
\Ss(\De) =\bigl\{v\in \Z^n\cap \De:   -v\in \De\bigr\} \less \{0\}
$$
of all integral symmetric points in $\De$.
The {\bf star-Ewald condition} for a vertex  $z$ says that
 there is a point $w\in \Ss(\De)$ 
that lies in precisely one of the facets through $z$ and is such that $-w$ 
lies on no facet through $z$.
As mentioned above, this condition is satisfied at $z$ exactly if all the points on the open line segment $C(z,0)$ from $z$ to $\{0\}$  can be displaced by probes.\footnote
{
To understand the conditions on $w$,  notice that 
the probes used to displace the points of the line  $C(z,0)$ have  base along the line $C(z,w)$
(which by hypothesis 
is contained in the interior of a facet $F_w$ through $z$) and direction $-w$.  The condition on 
$-w$ implies that the interior of the line $C(z,-w)$ lies in the {\it interior} of $\De$ so that the probes meet $C(z,0)$ {\it before} their halfway point.}
In particular, if $\De$ is a wedge with top and bottom facets $F_T, F_B$,
then to satisfy the  star-Ewald condition  at $z\in F_T\cap F_B$  the integer point $-w$ 
must lie on one of the other facets.  Because the union $F_T\cup F_B$ contains 
all the vertices of $\De$  and many
of its integer points, this condition is quite restrictive, and, as we now see, can fail to hold.

\begin{lemma}\label{le:sE} Let $\De$ be as in Equation \eqref{eq:sE}.
Then $\De$ does not satisfy the star-Ewald condition at the vertex $z= F_{123489}.$
\end{lemma}
\begin{proof}  Because the points ${\bf x} = (x_1,\dots,x_6)$ in $\De$ all satisfy the inequalities $x_i\ge -1$,  the coordinates of every point in  $\Ss(\De)$ line in the set  $\{0,\pm 1\}$.
Suppose that $w=(w_1,\dots,w_6)\in \Ss(\De)$
 lies in just one facet through $z$  while $-w$ lies in none of them.  Then
at most one of $w_1,\dots,w_4 $ is $-1$ and none is $1$.  If they are all $0$, 
then $w$ must lie on $F_8$ or $F_9$ so that precisely one of the equations
$ w_5+3w_6=1$ and  $w_5+2w_6=1$ holds.   Since 
$w_i\in \{0,\pm1\}$,  we must have $(w_5,w_6) = (-1,1)$.  
But then $w = (0,0,0,0,-1,1)$ does not lie in $\De$ because $w_5+3w_6>1$.
%But these equations have no integer solutions $w_i\in \{0,\pm1\}$ .
%this is impossible for integral $w_i$.    
 Therefore by symmetry we just need to consider the cases
 $$
{\rm (a)}\;\; \;w=(-1,0,0,0,w_5,w_6),\quad \mbox{ and }\quad {\rm (b)}\;\;\; w=(0,0,-1,0,w_5,w_6).
$$
In case (a),  since  $\pm w\in \De\less F_8$ we need
$
-1 +w_5+3w_6<1$ and $1-w_5-3w_6<1$.  This has the solution $(w_5,w_6) = (1,0)$. But then $w\in F_9$, which is not allowed.
A similar argument applies to case (b).
\end{proof}

  %%%%%%%%%%%%%%%%%%%%%%%%%%%%
  \subsection{Symplectic cutting}\labell{ss:cut}
 %%%%%%%%%%%%%%%%%%%%%%%%%%%%

Another useful way of constructing polytopes is by blow up.  As we show in more detail in \cite[\S3]{MT2}, blowing up along a face $f=F_I$ of codimension $k=|I|\ge 2$  adds a new face
$F_0$ to the polytope with conormal $\eta_0 = \sum_{i\in I}\eta_i$ 
and constant $\ka_0 = \sum_{i\in I}\ka_i - \eps$.  One can always do this for small $\eps>0$.  However, if $\De$ is monotone and one wants the blow up $\De'$ also to be monotone, then, because we need all the $\ka_j=1$, 
 one must take $\eps =k-1$.  In this case, the new facet $F_0$ is a $\De_{k-1}$-bundle  over $f$ whose fiber edges have affine length $k-1$, which is precisely the first Chern class of a line in the corresponding exceptional divisor.  In particular, there is a monotone blow up of  a vertex of a monotone polytope only if all edges through this point have affine length at least $n$.

In dimension $2$, blow ups have size $1$, and it is possible to make several such blowups on one polytope.  Indeed, one can blow up the triangle (the moment polytope of $\C P^2$) at all three of its vertices to obtain a monotone polytope.  Similarly,  in dimension $3$ one can make several disjoint monotone blowups provided they are along edges.  But
in dimension $3$ it is is not possible
to blow up two points simultaneously and in a monotone way.\footnote{
One can check this by listing the possibilities.  By \cite{B,WW}, there are only eighteen monotone polytopes in dimension $3$; they are all blow ups of polytopes obtained from simplices by forming suitable bundles.}
For example, the monotone $3$-simplex has edges of length $4$, while monotone blow ups (of vertices) have size $2$.   Thus if one did any two such blow ups one would create at least one singular (i.e. non simple) vertex.
\MS

\begin{quest}  Is there a monotone polytope   $\De$ of dimension $d>2$ for which  one can make at least two monotone and disjoint blow ups of points, or, more generally,  of any two faces of codimension $>2$?
\end{quest}
 
%	 \begin{rmk}\rm Every monotone toric symplectic  
%	 manifold is symplectically rationally connected, i.e. there is a nonzero genus $0$ Gromov--Witten invariant of the form 
%	 $\langle pt, pt,a_3,\dots,a_k\rangle_{0,\be,k}$, where $a_i\in H_*(M)$ and $\be\in H_2(M)$. If we can choose the class $\be$ so that $\om(\be)\le 2n-2$ then  one cannot embed two disjoint symplectic balls of capacity $n-1$, which implies one cannot make this toric blow up.  However,  this bound on $\om(\be)$ is not always satisfied; for example
%	 if $M=(S^2)^n$ is a product of $2$-spheres then the minimal class is $\be = \sum_i [(S^2)_i]$ with $\om(\be) = 2n$.  In this case, of course, the individual classes $\al_i: =[S^2]_i$ are uniruling classes, i.e. 
%	 an invariant  $\langle pt, a_2,a_3,\dots,a_k\rangle_{0,\al,k}$ is nonzero, which means that one cannot embed even one ball of size $n-1$.
%	 Hence one could make two monotone point blow ups 
%	only if the minimal uniruling class has size $\ge n$ 
%	and the minimal connecting 
%	class has size $\ge 2n-1$.
%	 \end{rmk}

In dimension $4$ is it not clear exactly what geometric constructions are needed to form all the monotone polytopes.  Obviously one can use bundles, or  wedges of lower dimensional (nonmonotone) polytopes.
Here is a monotone polytope formed by a different construction, that I again owe  to
Paffenholz \cite{Paf}.    

\begin{example}\rm  Let $\De$ be the $4$-dimensional cube 
$\{{\bf x}\in \R^4: -1\le x_i\le 1\}$.
Add the new facet $\sum x_i\ge -1$.  The new conormal $\eta_0 = (-1,\dots,-1)$ is parallel to the exceptional divisor that one would obtain by
blowing $\De$ up at its vertex $(-1,\dots,-1)$.  However, we take $\eps = 3$ to make a monotone blow up.  This means that we have cut out some of the vertices and edges of $\De$, 
though none of its facets.
The resulting polytope $\De'$ is smooth because none of the vertices of $\De$ lie on the new facet.  Really one should think of  $\De'$ not as a blow up but as the result of symplectic cutting; cf. Lerman \cite{L}.   As Paffenholz pointed out, $\De'$ is not a wedge because its vertices do not all lie on two facets, and it is not a bundle because it is not combinatorially equivalent to a product --- its facet $F_0$ has more vertices than any other.
\end{example}

Let us say that a  polytope is {\bf elementary} if removing any of its facets (i.e. deleting the corresponding inequality from the description given in equation \ref{eq:De} of $\De$) results either in a non simple or in an  unbounded polytope.
Clearly any polytope can be obtained from an elementary
 one by adding facets.  Adding a facet is the most general 
 possible cutting operation, where we no longer restrict the direction of the cut to $\sum_{i\in I} \eta_i$ for some face $F_I$.    If one wants to understand the  structure of monotone polytopes one might begin by asking about elementary ones.

 \begin{quest} What are the shapes of elementary monotone polytopes?  For example, is there any such polytope that is not a bundle or wedge?
\end{quest}

 %%%%%%%%%%%%%%%%%%%%%%%%%%%%%%%%%%%%%%%%%%%%%%%%%%%%%%%%
  \subsection{Bundles}\labell{ss:bun}
 %%%%%%%%%%%%%%%%%%%%%%%%%%%%%%%%%%%%%%%%%%%%%%%%%%%%%%%%

The general definition of bundle in the context of moment 
polytopes is rather complicated (see \cite[Def.~3.10]{MT1}), 
but that of bundle over the $k$-simplex
$\De_k$ is easy since the structure is determined by one \lq\lq slanted" facet $F_{N+k+1}$.   Note that in the following definition
we put the base coordinates last, since this is slightly more convenient and follows \cite{MT1}. 

\begin{defn}\labell{def:bun}   Write $\De$
 as in Equation \eqref{eq:poly},  and
normalize by assuming
  $\eta_i = -e_i, i=1,\dots,n$,
 where $e_i, i=1,\dots,n,$ forms the standard basis of $\ft=\R^n$.
  Then the  bundle $\De'$ with fiber $\De$ and base $\De_k$  
is determined (up to integral affine equivalence) by an 
integral $n$-vector $A : = (a_1,\dots,a_n)$ and constant 
$h: = \ka_{N+k}+\ka_{N+k+1}$ as follows. The polytope $\De'$   
  lies in $\R^{n+k} = \R^n\times \R^k$ and has conormals
\begin{gather}\labell{eq:bun}
  \eta_{i}' = (\eta_i,0,\dots,0), \;1\le i\le N,\;\; \;
  \eta_{N+i}' = -e_{n+i},\; i=1,\dots,k\;\;
  \\\notag
  \eta_{N+k+1}' = (a_1,\dots,a_n,1,\dots,1) = \sum_{i=1}^k\eta_{N+i} + A, 
\end{gather}
where $e_{n+i},  1\le i\le k,$ form the rest of  the standard basis  in $\R^{n+k}$. 
  The constants  $\ka_i, 1\le n\le N,$ are as before, we take $\ka_{N+i}=1, 1\le i\le k$, and choose  
 $\ka_{N+k+1}: = h-1$  large enough that the polytope is combinatorially a product of $\De$ with $\De_k$.
  \end{defn}
  
\begin{rmk}\labell{rmk:bun} \rm (i)
A bundle over $\De_1$   is formed from the product 
$\De\times \R$ by making two slices.  The bottom slice 
$F_{N+1}$ can be normalized to have equation $x_{n+1}= -1$, 
while the top $F_{N+2}$ is slanted by the vector $A$.
%; cf. Equation \eqref{eq:prod} below. 
If $M: = M_{\De}$, the
toric manifold
 $M_{\De'}$ is precisely $M_H$ where 
$H: =A=\sum_{i=1}^n a_i e_i$; see \cite[Ex.~5.3]{MT1}.
(As always, a formula like this involves some sign conventions; 
here we follow the sign choices in \cite{MT1}.)
\MS

\NI (ii) 
If the total space of a bundle (over an arbitrary base $\Hat\De$)
is monotone so are the fiber
and base. Even the meaning of the second part of this
statement needs clarification; see \cite[Lemma~5.2]{Mcpr}.
However, the first part is straightforward.  To prove it, 
notice that
 $\De\subset \R^n$ can be identified with the face 
 $$
 f': = 
 \De'\cap \bigl(\R^n\times \{(-1,\dots,-1)\}\bigr) = \bigcap_{1\le i\le k}F'_{N+i} 
 $$ 
 of $\De'$. Therefore the unit edge vectors from a vertex $V\in f'$ divide into two groups, the first group consisting of the edge vectors $e_{Vj}', j=1,\dots,n$  corresponding to the edge vectors $e_{Vj}$ in $\De$ and the second group given by $e_{Vj}': = e_{j}, j=n+1,\dots,n+k$.  Thus, using the notation of Lemma \ref{le:mono}, we have
 $$
 V+\sum_je_{Vj}' = u_0\in \R^n\Longleftrightarrow
V+\sum_{j\le n}e_{Vj} + \sum_i e_{n+i} = (u_0,0,\dots,0)\in R^{n+k}.
$$
Therefore the vertex-Fano condition at the vertex $V$ of $\De'$ readily implies this condition at the corresponding vertex $V$ of $\De$.  
\end{rmk}

\MS

\NI {\bf Proof of Proposition \ref{prop:monprod}.} 
We are given a monotone manifold of the form $(M_H,\om_H)$ where $H\in \ft$ generates a loop $\La_H$ of symplectomorphisms 
of the toric manifold $(M,\om,T)$.  To prove (i) we must show that if $H\ne 0$ then $H = \eta_j$, the conormal to one of facets of the moment polytope $\De$ of $M$.     
By Remark \ref{rmk:bun}(ii) the manifold $(M,\om,T)$ is
monotone. Moreover, if we identify $\De$ with the facet $F_{N+1} = \{x_{n+1}=-1\}$ 
of $\De'= \De_H$, each  vertex $V$ of $\De$ lies on a unique edge 
$\vareps_V$  that is parallel to $e_{n+1}$.  Choose $V_0\in \De$ 
such that this edge is shortest and then choose affine coordinates so that
 $V_0=(-1,\dots,-1)$, and the facets at $V_0$ have conormals $-e_i, i=1,\dots,n+1$.  Since the center point of $\De_H$ is $\{0\}$, all $\ka_j = 1$.

As in Remark \ref{rmk:bun}(i), the top facet of $\De_H$ is 
given by an equation of the form
$\sum_{i\le n} a_i x_i + x_{n+1} = \ka_{N+2}=1$, where $H = (a_1,\dots,a_n).$
 Therefore if $W = (w_1,\dots,w_n,-1)$ is a vertex in $\De \equiv F_{N+1}$, 
 the second endpoint of the edge $\vareps_W$ has 
 last coordinate equal to 
 $
1 -\sum_{i\le n} a_i w_i.$ Hence $\vareps_W$ has  length 
 $$
 \ell(\vareps_W) = 2-\sum_{i\le n} a_i w_i\ge \ell(\vareps_{V_0}) 
 =2+ \sum_{i\le n} a_i.
 $$
 But for each $i$ there is a vertex $W_i$ along the edge from $V_0$ in the direction of the $i$th coordinate axis. Thus $W_i=(-1,\dots,-1,w_i,-1,\dots,-1)$ where $w_i> -1$. Substituting $W=W_i$ in the above inequality, we find that $a_i\le 0$ for all $i$.

Now consider the vertex-Fano condition at the point 
$$
W = \vareps_{V_0}\cap F_{N+2} = (-1,\dots,-1,x_{n+1})
$$
on the top facet.
Because $\vareps_{V_0}$ is the shortest vertical edge, 
all edges through $W$ except for $-\vareps_{V_0}$ point in 
directions whose last coordinate is non-negative. 
(In fact, one can check as above that the unit vectors 
along these edges are  $e_i - a_ie_{n+1}$.) 
Hence, because 
$x_{n+1}\ge 0$,
we must have $x_{n+1}=0$ or $x_{n+1}=1$.
In the latter  case all these edges have zero last coordinate, 
which implies that $H=(a_1,\dots,a_n) = 0$.  Hence $F_{N+2}$ 
is parallel to $F_{N+1}$ and $\De$ is a product.
In the former case exactly one $a_i$ is nonzero, and the vertex-Fano condition $-1 - \sum a_i = 0$ shows that $a_i=-1$ 
Hence   $H = -e_i = \eta_i$ for some $i\le n$.  Thus $\La_H$ is the 
loop given by a rotation that fixes all the points in the facet $F_i$.  
This proves (i).

To prove (ii) we must show that the bundle formed from $H=\eta_i$
is never trivial.  Thus we must show  that such a  loop $\La_H$ 
is never contractible.  One way to prove this is to consider the Seidel 
 representation $\Ss$ of the group $\pi_1\bigl(\Ham(M,\om)\bigr)$ 
in the group $QH^{\times}$ 
of degree $2n$ units in the quantum homology ring of $(M,\om)$.
(For a definition of $\Ss$ in the toric context see \cite[\S2.3]{MT0}.)
 In the Fano case, it is easy to see that if $H = \eta_i$, we have $\Ss(\La_H)
= [F_i]\otimes \la$, where $\la$ is some unit in the Novikov 
coefficient ring of quantum homology, and $[F_i]$ denotes the 
homology class of $\Phi^{-1}(F_i)$, the maximal\footnote
{
i.e. the fixed point set on which the moment map $H\circ\Phi$ takes its maximum.}
 fixed 
point set of the loop $\La_H$; 
see for example \cite[Thm.~1.9]{MT0}.  
Since $\Ss(\la_H) \ne [M]$
 (the unit in $QH^{\times}$),
the loop $\La_H$ cannot be contractible.  

Similarly, if the loops $\La_{H_i}$ and $\La_{H_j}$ are homotopic
they must have equal images under $\Ss$ so that $[F_i]=[F_j]$.
But it is well known that the additive relations on $H_{2n-2}(M)$ have the form
$$
\sum_i \langle\eta_i,\xi\rangle [F_i] = 0, \qquad \xi\in \ft^*;
$$
see \cite{Tim} for example.  Hence $[F_i]=[F_j]$ iff these two 
facets are equivalent in the sense used here.  \QED

We end with a question.

\begin{quest} Is there a monotone polytope $\De$ that supports an essential mass linear function?
\end{quest}

Note that our constructions for monotone polytopes tend to destroy
essential mass linear functions.  For example, if $H$ is an essential mass linear
function on $\De$ and $\De'$ is the wedge of $\De$ along some facet then
$H$ does not in general induce an  essential mass linear function on $\De'$.
 A similar statement is true for bundles; if $\De'\to \Hat\De$ is a bundle with fiber $\De$, then essential mass linear functions on $\De$ do not usually extend to mass linear functions on $\De'$:  explicit examples are given in  \cite[\S3]{MT2}.

%%%%%%%%%%%%%%%%%%%%%%%%%%%%
  \subsection{An example of uniqueness}\labell{ss:uni}
 %%%%%%%%%%%%%%%%%%%%%%%%%%%%

We now prove
Proposition \ref{prop:uni}.  This   states that there is a
 unique toric structure on the product $(M,\om): = (\C P^k\times \CP^m,\om_k\oplus \la \om_m)$, if $k\ge m\ge2$, or if $k>m=1$ and $\la\le 1$, or if $k=m=1$ and $\la =1$.  Notice that all monotone products of projective spaces satisfy these conditions.
 (Recall that we have normalized $\om_k$ so that its integral over a line is $1$.)

Suppose that $\De$ is the moment polytope for some toric structure on $M=\C P^k\times \CP^m$.  Then $\De$ has precisely $k+m+2 =
\dim \De + {\rm rank\,} H^2(M)$ facets.
Hence by  Timorin \cite[Prop.~1.1.1]{Tim}, $\De$ is combinatorially equivalent to a product 
of two simplices.  Therefore, because $\De$ is smooth,
 \cite[Lemma~4.10]{MT1} implies that $\De$ is a   $\De_r$-bundle over 
$\De_s $ for some $r,s$. 
Therefore $M$ is a $\C P^r$-bundle over $\C P^s$.  Hence $H^2(M;\Z)$ 
contains an element $\al$ such that $\al^{s+1}=0$ while $\al^s\ne 0$.  It follows that 
$s=k$ or $s=m$.  

Let us now suppose that $\De$ is not the trivial bundle, i.e.
 some $a_i\ne 0$ in the presentation
described in Definition \ref{def:bun} for a bundle over $\De_s$.
For each vertex $V$ of the fiber $\De_r$, there is an $s$-dimensional face $f_V$  of $\De$ that is affine equivalent to  $V\times \mu_V\De_s$ for some scaling constant 
$\mu_V$.  Choose $V$
so that $\mu_V$ % (resp. $\mu_{v'}$ 
is  minimal.  
As at the beginning of \S\ref{s:mono}, choose coordinates on $\R^{r+s}$ so that $V = (-1,\dots,-1)$ and so that the edges from $V$ point in
the directions of the coordinate axes.   Then, as in the proof of Proposition \ref{prop:monprod}, each $a_i\le 0$.

Now let us calculate $H^*(M;\Z)$ using the Stanley--Reisner presentation  
\begin{equation}\labell{eq:SR}
 \Z[x_1,\dots,x_{r+s+2}]/\bigl(P(\De) + S(\De)\bigr),
\end{equation}
 where the additive relations
 $P(\De)$ are $\sum_i\langle\eta_{i},e_j\rangle x_i = 0$ (where $e_1,\dots,e_n$ is a basis for $\ft^*$),
  and the set $S(\De) $ of multiplicative relations is
  $\prod_{i\in I}x_i=0$, where $I$ ranges over all  minimal subsets 
 $I\subset \{1,\dots,N\}$ such that the intersection $F_I: = \cap _{i\in I} F_i$ is empty.  Thus in the case at hand there are two 
 multiplicative relations,
one from the conormals $(\eta_i,0)$, where $1\le i\le N = r+1$, and the other from the conormals  $\eta_{N+i}, 1\le i\le s+1$.
 The  relations for the  facets $F_{N+i}, 1\le i\le s+1,$
 show that $x_{r+2} = \dots = x_{r+s+1} = \al,$  say, and $\al^{s+1}=0$.  Similarly for the  facets $F_i, 1\le i\le r+1$, we find
 $$
 -x_i + x_{r+1} + a_i\al=0, \;\; i=1,\dots,r,\quad  \prod_{i=1}^{r+1}  x_i=0.
 $$
 Thus, if we write $x_{r+1}: = \be$ and  define $a_{r+1}:=0$ , we find  
\begin{equation}\labell{eq:berel}
 0 = \prod_{i=1}^{r+1}(\be+a_i\al) = 
 \be^{r+1} +\si_1\be^r \al + \dots + \si_r \be \al^{r},
\end{equation}
where $\si_1: = \sum_i a_i$, and, more generally, $\si_k$ is the value of the $k$th elementary symmetric polynomial on $(a_1,\dots,a_r,0)$.
 
By assumption, there are generators  $\al_0,\be_0\in H^2(M;\Z)$ 
so that
$\al_0^{s+1} = 0 = \be_0^{r+1}$.  Therefore, for some 
$A,B,C,D\in \Z$ with $ AD-BC=1$, we must have
$$
( A\al+B\be)^{s+1} = 0 = (C\al+D\be)^{r+1}.
$$
%Further, 
% there is a ring isomorphism $H^*(M;\Z)\to H^*(M;\Z)$,
% that must have the form
%$$
%\al_0\mapsto A\al+B\be,\quad \be_0\mapsto C\al+D\be,
%$$
% where $A,B,C,D\in \Z,$ and $ AD-BC=1$.
 We now divide into cases, and show in each case
 if some $a_i$ is nonzero 
then the conditions in Proposition \ref{prop:uni} must hold.
\MS

\NI {\bf Case 1:}  $1<s<r$. 

In this case, the ring $H^*(M)$ is freely generated by $\al,\be$ in degrees $\le 2s$ and there are two relations  of degree $2s+2$, namely 
$$
\al^{s+1}=0,\quad (A\al + B\be)^{s+1} = 0.
$$
If $B\ne0$, these relations are different so that $H^{2s+2}(M)$ has rank 
$s$ instead of $s+1$.  Therefore $B=0$ and $A=\pm 1$.   By changing the sign of $\al_0$ we may suppose that $A=1$, so that
$D=1$.    Then we have $(C\al + \be)^{r+1} = 0$.  
Again, this must agree term by term with equation \eqref{eq:berel}, 
once we substitute $\al^{s+1}=0$.
Equating coefficients for $\be^{r-i+1} \al^i$ with $i=1,2$, we need
$
\sum_{i\le r} a_i = (r+1)C,$ and $\sum _{i\ne j} a_ia_j =  (r+1)rC^2$.  Hence we also need 
$\sum a_i^2 = (r+1)C^2$.
But the last inequality together with Cauchy--Schwartz gives
$$
|\sum a_i|\le \sqrt r \sqrt{\sum_ia_i^2}\le \sqrt r\sqrt{r+1} |C|,
$$
a contradiction. Therefore this case does not  occur.
\MS

\NI {\bf Case 2:}  $s>r$.  

In this case, the relation
$(C\al+D\be)^{r+1}$ must be a nonzero multiple of equation \eqref{eq:berel}.  But, if $C\ne 0$, the coefficient of $\al^{r+1}$ is nonzero in the first equation, while it vanishes in \eqref{eq:berel}.  Therefore
$C=0$, so that $\si_1 = \sum a_i = 0$.  But each $a_i\le 0$ by construction.
Hence we must have  $a_i=0$ for all $i$.  Hence, again this case does not occur.\MS

\NI {\bf Case 3:}  $s=r>1$

In this case we have four relations in degree $2r+2\ge 6$, namely  equation \eqref{eq:berel} and
$$
\al^{r+1} = 0, \quad  (A\al+B\be)^{r+1}=0,\quad 
(C\al+D\be)^{r+1}=0
$$
that must  impose just two linearly independent conditions.   Since $\al^{r+1} = 0$ is independent from \eqref{eq:berel}, the other two equations must be combinations of these.  By permuting $\al_0,\be_0$ if necessary,
 we can suppose that 
$A\ne 0, D\ne 0$.    But then if we put $\al^{r+1} = 0$ in the relation
$(C\al+D\be)^{r+1}=0$, we must get $D^{r+1}$ times the equation \eqref{eq:berel}.
Comparing coefficients of $\be^r\al$ and $\be^{r-1}\al^2$ we find
$$
(r+1) C = D\si_1, \quad \tfrac {r(r+1)}2C^2 = D^{2}\si_2,
$$
which, as in Case 1,  is impossible unless all $a_i=0$.   Therefore this case also does not occur.\MS

\NI {\bf Case 4:} $r\ge s=1$.

% If $M$ is a product there must be an element $z=y+kx\in H^2(M;\Z)$ 
% such that $z^{r+1} = 0$.  
%Let us now suppose that $r< s$ so that $x^{r+1}\ne 0$.  We then find  
% $$
% \sum a_i = -k(r+1), \;\;\sum_{i\ne j} a_ia_j = k^2(r+1)r/2,\;\;\dots, \;\;k^{r+1} = \pm 1.
% $$
%Since  $a_{r+1} = 0$ by assumption, there are no nonzero solutions to these equations with $a_i\ge 0$ for all $i$.  Hence this case does not occur.
%
%So suppose that $r\ge s$.   
%Since each $a_i\ge 0$ and $-\sum a_i = k(r+1)$, we must have
%$\sum a_i \ge r+1$.    
Let us go back to the polytope $\De$ and look 
at the  face $f_V\cong \mu_V\De^s$ at our chosen vertex $V$.
 Every edge $\vareps$ in $f_V$ has first Chern class given by\footnote
 {
 Here $c_1(\vareps)$ is more correctly described as the first Chern class of the restriction of the tangent bundle $TM$ to the $2$-sphere $\Phi^{-1}(\vareps)$.  The paper \cite{KKP} describes how to  calculate $c_1(\vareps)$ when $n=2$. See also 
  Remark \ref{rmk:finit} below.}
\begin{equation}\labell{eq:c}
c_1(\vareps) = s+1+\sum a_i.
\end{equation}
 Now observe that the submanifold $\Phi^{-1}(f_V)$ is a section of the bundle
$M\to \C P^s$, so that the 
$2$-sphere $\Phi^{-1}(\vareps)$ lies in a homology class of the form 
$q L_r + L_s$, where $L_i$ denotes the line in $\C P^i$.

Now observe that if $r>s=1$ we must have $r=k$ and $s=m$, while if 
$r=s=1$ we may  assume that $r=k$ and $s=m$.  Then, in both cases,
we have $\om(L_r) = 1$ and $\om(L_s) = \la$.    Hence $
\om(q L_r + L_s)  = q + \la> 0$. 
On the other hand if some $a_i<0$ then  $c_1(\vareps) < s+1$ so that $q<0$.  
Therefore  this case does not occur when $\la\le1$.

Finally note that if $r=s=1$ we can interchange the roles of $r$ and $s$, replacing $\la $ by $1/\la$.
Therefore, when $k=m=1$ our argument rules out the existence of nontrivial bundles 
only in the case $\la=1$.

This completes the proof.

\begin{rmk}\labell{rmk:trivk}\rm  As is clear from Definition \ref{def:bun}, 
toric bundles over $\De_k$ and with fiber $\Tilde\De$ of dimension $r$
are determined by one vector $H = -(a_1,\dots, a_r)$ that generates a 
circle action $\La_H$ on the fiber $(\Tilde M,
\om): = (M_{\Tilde \De},\om_\ka)$.  
It is tempting to think that this bundle is trivial as long as this 
circle contracts in $\Ham({\Tilde M},\om)$. But as we saw above, 
this clearly need not be so when $k>1$.  
For example,  if $\Tilde\De = \De_r$ then $H=(1,0,\dots,0)$  
generates a circle $\La_H$ that lies in $SU(r+1)$.    Since
$\pi_1(SU(r+1)) \cong \Z/(r+1)\Z$, we find that    
 $\La_{(r+1)H}$ contracts in $SU(r+1)$ and hence in $\Ham (\C P^r)$.
On the other hand, by Proposition \ref{prop:uni} the bundle is nontrivial when $k>1$.
 
To understand this notice that, if $\La_H$ contracts, then
the classifying map $\C P^k\to B\Ham(\Tilde M,\om))$ induces the null map on the $2$-skeleton $\C P^1\subset \C P^k$. When we contract this $2$-sphere, we 
get further obstructions to the null homotopy of the whole map.  These 
obstructions are explained in K\c edra--McDuff; see \cite[Thm~1.1]{KM}. 
We show there that 
the existence of  the contractible
 circle $\La_H$  in $\Ham(M,\om)$ creates a  nonzero element\footnote
 {
 It is detected by  a characteristic class very similar to those used by Shelukhin; cf. equation \eqref{eq:I0} below.}
 (a kind of Samelson product) in  
 $\pi_3\bigl(\Ham(M,\om)\bigr)\cong \pi_4\bigl(B\Ham(M,\om)\bigr)$,  that has nonzero pullback
 under the classifying map $\C P^k\to B \Ham(M,\om)$ of this bundle. 
  Thus the bundle is nontrivial.  This makes it unlikely that the total space
could ever be diffeomorphic to a product, though it 
does not completely rule it out without further argument.\end{rmk}

  %%%%%%%%%%%%%%%%%%%%%%%%%%%%
  \section{Questions concerning finiteness}\labell{s:fin}
 %%%%%%%%%%%%%%%%%%%%%%%%%%%%

 %%%%%%%%%%%%%%%%%%%%%%%%%%%%
  \subsection{Finite number of toric structures}\labell{ss:fin}
 %%%%%%%%%%%%%%%%%%%%%%%%%%%%

 \begin{prop}\labell{prop:fin}  Let $(M,\om)$ be  a $2n$-dimensional   symplectic manifold.  Then, the number of distinct toric structures on $(M,\om)$ is finite, 
 where we identify 
  equivariantly symplectomorphic actions.
  \end{prop}
  
  \begin{proof}   Let $\Phi:M\to \R^n$ be the
   moment map of some toric structure on 
  $(M,\om)$ with image $\Phi(M)=:\De$.     The number $N$ of facets of
  the polytope $\De$ is $n+\dim H^2(M;\R)$.
 We first show  that $\De$ is determined by 
 the classes $x_i\in H^2(M;\Z), i=1,\dots,N,$
  that are Poincar\'e dual to the divisors 
  $\Phi^{-1}(F_i)$ corresponding to the facets $F_i$.  Then we will show that   
  these classes $x_i$ lie in a finite subset of $H^2(M;\Z)$.
  
  To prove the first statement,
  number the $x_i$ so that $x_1x_2\dots x_n\ne 0$  and  $e_1: = -x_1,\dots,e_n: = -x_n$ form a basis for $H^2(M;\Z)$.
  Then the Stanley--Reisner presentation of $H^*(M)$ (cf. equation \eqref{eq:SR})  implies that
    the coordinates of the conormals for the other facets can be read off from the linear relations between   the $x_i, i=1,\dots,N$.   (Recall that we always assume that the conormals are primitive integral vectors, i.e. that their 
    coefficients have no common factor.)
    
    Therefore it remains to determine the support constants $\ka_i$.   Because of the translational invariance of $\De$, the  first $n$ of these can be chosen at will.  Once these are  chosen, the other $\ka_i$ can be determined
     by looking at  a suitably ordered set of edge lengths.  
 To see this,  let us set $\ka_i = 0, i\le n,$ so that 
    $$
    v_0:=\cap_{i=1}^n F_i =(0,\dots,0).
    $$
Suppose that $v_1$ is connected to $v_0$
 by the edge $\vareps_j$ that is transverse to 
 $F_j$ at $v_0$ for some $j\le n$.  Then 
$\vareps_j = \cap_{i\in I} F_i$ where $I: = \{i\le n, i\ne j\}$,
 and its affine length is
    $$
    \ell(\vareps_j) = \int_{\vareps_j} [\om] = \int_M x_I\; [\om],\quad\mbox{where } x_I: = \prod_{i\in I} x_i.
     $$
If the other endpoint of $\vareps_j$ is transverse to $F_k$, then $\ka_k$ is determined by $\ell(\vareps_j)$.  Proceeding in this way, we can find $\ka_k$ first for all facets joined to $v_0$ by one edge, then for those joined to $v_0$ by a path consisting of two edges, and so on.

Therefore it suffices to 
 show that there are a finite number of possibilities for  these classes $x_i\in H^2(M;\Z), i=1,\dots,N$.
  Following a suggestion of Borisov,\footnote{Private communication.} let us look at the Hodge--Riemann form  on $H^2(M;\R)$ given by
  $$
  \langle \al,\be\rangle: = \int_M \al\, \be \,\om^{n-2}.
  $$
  By the Hodge index theorem, this is nondegenerate  of type
 $(1,-1,\dots, -1)$; in other words it is negative definite on the orthogonal complement to $[\om]$.  (A nonanalytic proof of this result for  toric manifolds may be found in Timorin \cite{Tim}.)  
 Write $x_i = y_i + r_i[\om]$ where $  \langle y_i,\om\rangle = 0$ and $r_i\in \R$.   Then each $r_i>0$, since   
  $$
  r_i  \langle \om,\om\rangle =  \langle x_i,\om\rangle =
  \int_{F_i} \om^{n-1} >0
  $$
  because it is a positive multiple of the $\om$-volume of the 
  K\"ahler submanifold $\Phi^{-1}(F_i)$.
  Further, because   $c_1(M) = \sum x_i$ by Davis--Januszkiewicz
  \cite{DJ}, we have
  $$
  \sum_i r_i\, \langle \om,\om\rangle = \sum_i   \int_M x_i\, \om^{n-1} =
  \int_M c_1(M) \,\om^{n-1} = : C \langle \om,\om\rangle.
  $$
 Therefore, each $r_i< C$, 
 so that $\sum r_i^2< NC^2$.
  Finally, because $c_1^2 - 2c_2 = \sum x_i^2$ we have
  $$
A:= \int_M (c_1^2 - 2c_2)\,\om^{n-2} =
  \sum_i   \langle x_i,x_i\rangle= \sum r_i^2 + \sum_i\langle y_i,y_i\rangle.
  $$
  Since each $\langle y_i,y_i\rangle\le 0$ by the Hodge index theorem,
  we find that
  $$
  0\le -\sum_i\langle y_i,y_i\rangle \le NC^2/V^2-A.
  $$
  Therefore the integral classes $x_i$ lie in a bounded subset of 
  $H^2(M;\R)$.  
  Thus they are all contained in a finite subset of $H^2(M;\Z)$.
\end{proof}

\begin{rmk} \labell{rmk:finit} \rm There are various  elementary proofs of  
finiteness when $[\om]$ is integral.   Perhaps the simplest is again due to Borisov, who pointed out the the following argument.  Normalize $\De$ so that one vertex is at the origin and the edges from it point along the positive coordinate axes.   Denote by  $S_i$ the $(n-1)$ simplex in the hyperplane $\xi_i = 0$ with edges of unit length, and suppose that $v = (a_1,\dots,a_n)\in \Z^n$ is some vertex   of $\De$.  Then the volume of the cone spanned by $S_i$ and $v$ is $a_i/n!$.  Since this cone lies in $\De$,
 the coordinates of $v$ are bounded by the volume $V$ of $\De$.
Therefore, the vertices  lie in a bounded subset of the lattice $\Z^n$ whose
 size is determined by $V = \frac 1{n!} \int_M\om^n$.

Another approach is first to note that the number and affine lengths of the edges are bounded by some constants $K, L$  because  
the sum of their lengths is $\int_M c_{n-1}\om$ and each edge has length at least $1$.    It then follows that the geometry of each edge $\vareps$ is bounded.  To see this, note that this geometry is  determined by the 
Chern numbers $c_{F_i}(\vareps)$ of the normal line bundle to $F_i$ along 
$\Phi^{-1}(\vareps)$, where $F_1,\dots,F_{n-1}$ are the facets containing $\vareps$. Because each edge has length between $1$ and $L$ and each $2$-face is a convex polygon, we must have $c_{F_i}(\vareps)\le L$ for each such $i$.  
But  $\int_{\Phi^{-1}(\vareps)} c_1(M) = 2 + \sum_i c_{F_i}(\vareps)$.
It follows that the $c_{F_i}(\vareps)$ are bounded above and below.
 Since $\De$ is made by putting together at most $K$ edges, there are again only finitely many possibilities for $\De$.
\end{rmk}

\NI {\bf Proof of Theorem \ref{thm:fin}.}
The proof of Proposition \ref{prop:fin} used only cohomological facts about $M$.  The number of facets of $\De$ is determined by the rank of 
$H^2(M)$.  We also needed to know $\int_M\om^n$, $\int_M c_{1} \om^{n-1}$ and $\int_M (c_1^2-2c_2)\,\om^{n-2}$.  But, once one knows the classes $c_1,c_2$ and $[\om]$, these integrals are determined by  
the integral cohomology ring.   This holds because
 there is a unique generator $u$ of $H^{2n}(M;\Z)$ such that $\om^n = \la u$ for some $\la>0$, and then
an integral such as $\int_M c_{1} \om^{n-1}$ is equal to  $a\in \R$, 
where $
c_{1} \om^{n-1} = a \,u$.
This completes the proof.\QED

  %%%%%%%%%%%%%%%%%%%%%%%%%%%%
  \subsection{Manifolds with more than one toric structure: blow ups}\labell{ss:many}
 %%%%%%%%%%%%%%%%%%%%%%%%%%%%

One easy to way to construct different toric structures on a symplectic manifold $(M,\om)$   is by blowing up.   Suppose given a toric structure on $(M,\om)$ with moment map $\Phi:M\to\De$.   
As we show in more detail in \cite[\S3]{MT2}, blowing up along a face $f=f_I$ of codimension $k=|I|\ge 2$  adds a new facet
$F_0$ to the polytope with conormal $\eta_0 = \sum_i\eta_i$ 
and constant $\ka_0 = \sum_i\ka_i - \eps$.  The new moment polytope
$\De_f$ is $\De\less Y_{f,\eps}$, where 
$$
Y_{f,\eps}=\bigl\{\xi\in \De: \langle \eta_0,\xi\rangle > \ka_0\bigr\}.
$$
This is a smooth moment polytope for small $\eps>0$.  The corresponding symplectic manifold $(\Tilde M_f, \Tilde \om_\eps)$  is formed from $(M,\om)$ by excising $\Phi^{-1}(Y_{f,\eps})$ and collapsing the boundary along its characteristic flow.  This is an example of symplectic cutting; cf Lerman \cite{L}.  
If $f=v$ is a vertex,  we call the resulting toric manifold a {\bf one point toric blow up of weight} $\eps$.    The underlying symplectic manifold is called the\footnote
{
A subtle point is concealed here. For most manifolds it is not known whether there is $\eps_0>0$ such that the space of symplectic embeddings of a ball of size $\eps\le \eps_0$ into $(M,\om)$  is connected.  Since each such embedding gives rise to a symplectic blow up (see \cite{McP}), it is not known whether all  sufficiently small one point blow ups are symplectomorphic.  In the toric case,  this problem does not arise
since we have given a unique way to do such a blow up at each vertex. 
Even if one allows an ostensibly more general process by using equivariants embeddings, the invariance of the image (see  Pelayo 
\cite[Lemma~2.1]{Pel}) 
shows that this gives nothing new.}
 one point blow up of $(M,\om)$.

 \begin{lemma} \labell{le:torbl}
  Let $(M,\om,T)$ be a toric manifold with moment polytope $\De$.  Then there is $\eps_0>0$ such that if $0<\eps<\eps_0$ all of its one point toric blow ups  
  of weight $\eps$ are symplectomorphic.
 \end{lemma}
 \begin{proof}[Sketch of proof]
 In this case 
$\Phi^{-1}(Y_{f,\eps})$ is the image of a standard ball $B^{2n}(\eps)$ of radius $\sqrt {\eps/\pi}$.
 Because the Hamiltonian group acts transitively on $M$, it is easy to see that given any two  symplectic embeddings  
 $B^{2n}(\eps')\to M$ one can find $\eps_0\in (0,\eps')$ such that their restrictions to $B^{2n}(\eps_0)$ --- and hence to any smaller ball --- are isotopic.  This implies that the corresponding blow up manifolds 
 are symplectomorphic; see for example \cite{McP}.  To complete the proof, it remains to observe that there are a finite number of toric blow ups.
\end{proof}

Similarly, if one blows up along faces $f,f'$ for which the inverse images
$\Phi^{-1}(f)$ and $\Phi^{-1}(f')$ are Hamiltonian isotopic, the resulting 
blow ups are symplectomorphic for small enough $\eps$.

\begin{rmk}\rm  There are many interesting 
questions here about exactly how big one can take $\eps_0$ to be; 
cf. the discussion in Pelayo \cite[\S3]{Pel}.
\end{rmk}

\begin{defn}\labell{def:equiverr} 
Two vertices of $\De$ are said to be {\bf equivalent} if there is 
an integral affine self-map of $\De$ taking one to the other. Further we define $\Aff(\De) : = \Aff(\De(\ka))$ to be the group of all integral affine self-maps of $\De(\ka)$, and $\Aff_0(\De)$ to be the subgroup that is generated by reflections that interchange equivalent facets.
\end{defn}

As explained in the discussion after  Definition \ref{def:equiv}, the elements of $\Aff_0(\De)$ lift to elements 
in the Hamiltonian group of $(M,\om)$ while the elements in $\Aff(\De)\less \Aff_0(\De)$  lift to sympectomorphisms that are not isotopic to the identity; in fact, because they interchange nonequivalent facets, they act nontrivially on $H^2(M)$.  (See also \cite{Mas}.)   
Observe that $\Aff(\De)$ depends on $\ka$, while $\Aff_0(\De)$ does not.  (In the terminology of \cite{MT1}, $\Aff_0(\De)$ consists of robust transformations.)     In particular, the question of which vertices  of $\De(\ka)$ are equivalent depends on $\ka$.     Explicit examples of this are provided by blow ups of $\C P^2$. Note also that, 
when $\ka$ is generic, the two groups $\Aff(\De(\ka))$ and $\Aff_0(\De(\ka))$ coincide because different homology classes in 
$H_{2n-2}(M)$ are distinguished by the $\om$-volume of their representatives.

The next corollary follows by combining these remarks with Lemma \ref{le:torbl}.

\begin{cor}\labell{cor:torbl}  Let $(M,\om,T)$ be a toric manifold whose moment polytope has $k$ pairwise nonequivalent vertices. Then there is $\eps_0>0$ such that 
for all $0<\eps<\eps_0$ the one point $\eps$-blow up of $(M,\om)$ has at least $k$ different toric structures.   
\end{cor}
\begin{proof}  Let $\De$ be the moment polytope of $(M,\om,T)$, and let $\De'_\eps$ and $\De''_\eps$ be the polytopes obtained by blowing up $\De$ at two nonequivalent vertices
$v'$ and $v''$ by some amount $\eps>0$.     Delzant's theorem \cite{Del}
 states that  a toric manifold is determined up to equivariant symplectomorphism
 by the integral affine equivalence class of its moment polytope.  
Therefore if the toric manifolds corresponding to these two blow ups are 
equivariantly symplectomorphic there is an integral affine transformation $A_\eps$ taking $\De'_\eps$ to $\De''_\eps$.
We may suppose that 
 $\eps$ is less than half  the length of the shortest edge of $\De$.  Then the only facet of $\De'_\eps$ with all edges of length $\eps$ its the exceptional divisor $F_0'$.  Since a similar statement holds for $\De''_\eps$, the facet $A_\eps(F_0')$ must be the exceptional divisor $F_0''$ of $\De''_\eps$.  But the quantities 
 $\eta_i$ and $\ka_i$ that determine
 the other facets of $\De'_\eps$ and $\De''_\eps$ are independent of $\eps$ and must be permuted by $A_\eps$.  
Hence $A_\eps$  induces a self-map of $\De$ which is independent of $\eps$ and takes  $v'$ to $v''$.  Thus $v'$ is equivalent to $v''$, contrary to the hypothesis.  
Thus the toric manifolds obtained by blowing up two nonequivalent vertices are not 
  equivariantly symplectomorphic.  One the other hand, if $\eps>0$ is sufficiently small the underlying manifolds are symplectomorphic by Lemma \ref{le:torbl}. 
\end{proof}

We begin the proof of Proposition~\ref{prop:torbl}
by considering the following special case.  

\begin{lemma}\labell{le:equivf0}  Suppose that $\De$ is a nontrivial and generic 
bundle over $\De_1$.
 Then $\De$ has at least two 
inequivalent vertices.
\end{lemma}
\begin{proof}   By Definition \ref{def:bun} the structure of $\De$ is determined by the
vector $A = -(a_1,\dots,a_r)$ and the length $h$.  The bottom and top facets $F_{N+1}$ and $F_{N+2}$ are equivalent, and, as explained in the proof of Proposition \ref{prop:monprod}, 
 the numbers $h-a_1,\dots h-a_r, h-a_{r+1}$ (where $a_{r+1}: = 0$ as in \S\ref{ss:uni}) are the lengths of the vertical edges of $\De$, i.e. those that are parallel to $e_{r+1}$, going between   the bottom and top facets.  
  Since $A\ne 0$, there are two vertices $v_1,v_2$ on the bottom facet $F_{N+1}$ that are the endpoints of vertical edges of different lengths.   Since $\De$ is generic, we can choose $h$ so that every affine transformation of $\De$ must preserve the set of vertical edges (possibly changing their orientation), because there are no other edges of precisely these lengths.  Hence $v_1$ and $v_2$ 
  cannot be equivalent.  
\end{proof}

\begin{rmk}\rm 
It is possible that there are just two equivalence classes of vertices.  For instance $\De$ might be a $\De_2$-bundle over $\De_1$ with $A = -(0,1)$.
\end{rmk}

%The proof of  Proposition~\ref{prop:torbl} is based on 
The  following lemma
generalizes  Theorem 1.20 of \cite{MT1} which imposes the 
extra condition that  $F_I: = \cap_{i\in I}F_i=\emptyset$ for all equivalence classes $I$ and concludes that $\De$ is a product of simplices.

\begin{lemma}\labell{le:equivf}   Suppose that each facet $F_i$ of $\De$ is  equivalent to some other facet $F_j$.  Then, if $\ka$ is generic, either  $\De(\ka)$ is a product of simplices or $\De(\ka)$ has at least two nonequivalent vertices.
\end{lemma}
\begin{proof} If $\De$ has   dimension $2$,  then it must either be $\De_2$ or $\De_1\times \De_1$.
Now assume inductively 
that the lemma holds for all polytopes of dimension $\le n-1$, where $n: = \dim \De$.  

Suppose first that there is some equivalence class $I$ with $|I|\ge 3$ and renumber the facets so that  $\{1,2\}\subset I.$  Then by Proposition 3.17 of \cite{MT1}, $\De$ is the $1$-fold expansion of the facet $F_2$ along its facet $F_{12}$.  In particular $F_2$ is pervasive, i.e. meets all other facets.   Lemma~3.27 in \cite{MT1} states that, when $F_2$ is pervasive,
 two facets $F_j,F_k$, where $j,k\ne 2$,  are equivalent in $\De$ exactly if the facets $F_{2j}, F_{2k}$ are equivalent in $F_2$.  
Because $|I|\ge 3$, this implies that there is $j\in I\less \{1,2\}$
such that $F_{21}\sim F_{2j}$.   Hence,  the equivalence classes of facets of 
$F_2$ all have more than one element. Therefore, by the inductive hypothesis $F_2$  is either a product of simplices  or has at least two nonequivalent vertices.

In the former case,  $\De$ is the expansion of $F_2=\De_{k_1}\times \dots\times \De_{k_p}$ along a facet $F_{12}$ that we may assume to have the form
$F\times \De_{k_2}\times \dots\times \De_{k_p}$ for some facet $F$ of $\De_{k_1}$.
It is now easy to check from  the definition of expansion that 
$$
\De=
\De_{k_1+1}\times \dots\times \De_{k_p}.
$$
 
 In the latter case,  there are at least two vertices $v_1,v_2$ of $F_2$ that are not equivalent under $\Aff(F_2) = \Aff_0(F_2)$.  
 It suffices to show that they are not equivalent under $\Aff_0(\De)$.  
 Suppose not, and let $\phi\in \Aff_0(\De)$ be such that $\phi(v_1) = v_2.$  Then $\phi(F_2)\ne F_2$.  
 But because $\phi\in \Aff_0(\De)$  we must have $\phi(F_2)\sim F_2$.  Let $\al\in \Aff_0(\De)$ be the reflection that interchanges the facets $\phi(F_2)$ and $F_2$.  Then $\al\circ\phi(F_2)=F_2$.   
 But $v_2\in F_2\cap \phi(F_2)$ is fixed by $\al$.  Hence 
 $v_1$ and $v_2$ are equivalent in $F_2$, contrary to hypothesis.  
 This completes the proof when there is some equivalence class with $>2$ facets.

It remains to consider the case when all equivalence classes have two elements. 
Suppose there is such an equivalence class 
 $I = \{1,2\}$ with $F_{12}=\emptyset$.
 Then again each equivalence class of facets of  $F_2$ has at least $2$ elements, and  
 \cite[Prop.~3.17]{MT1} implies that  $\De$ is an $F_2$-bundle over $\De_1$.  If 
 this bundle is nontrivial, then Lemma \ref{le:equivf0} implies that $\De$ has at least $2$ nonequivalent vertices.  
 If it is trivial, then either 
  $F_2$ (and hence also $\De$) is a
 product of simplices,  or we can use
 %the we can apply the inductive hypothesis to $F_2$ to find 
 the two nonequivalent vertices of $F_2$ supplied by the 
 inductive hypothesis to find two such vertices of $\De$.
  
 The remaining  possibility is that 
each  equivalence class consists of precisely two
 intersecting facets $F_i, F_i', 1\le i\le \ell $.  In this case,  the proof is completed by Lemma \ref{le:pair} below.
 \end{proof}

 \begin{lemma}\labell{le:pair}  Suppose that $\De$ is a polytope such that each  facet is equivalent to at most  one other.  Suppose further that each pair of equivalent facets intersects.  Then, for generic $\ka$, the polytope $\De(\ka)$ has at least five
  nonequivalent vertices.
 \end{lemma}
 \begin{proof}  
 Pick one facet $F_i, 1\le i\le \ell, $ from each equivalence class 
 with more than one element,
 and denote the other facets in these equivalence classes by $F_i'$ where $F_i\sim F_i'$.   
We first claim that  the face   
 $f: = F_1\cap\dots\cap F_\ell$ has the following properties:
 \begin{itemize}
 \item[(a)] $f\ne \emptyset$;
 \item[(b)] $\De$ is made from $f$   by expanding once along each of the facets
 $F_i'\cap f$;
 \item[(c)] no two facets of $f$ are equivalent;
 \item [(d)]  $f$ has at least $5$ vertices.
 \end{itemize}
  We prove this by induction on $\ell$.  
 If $\ell = 1$,  (a) is clear and (b)  holds by 
  \cite[Prop.~3.17]{MT1} (which states that when $F_1\cap F_1'\ne \emptyset$, 
  the polytope
  $\De$ is the expansion of  $F_1$ along $F_1\cap F_1'$).  Therefore $F_1$ is pervasive, so that we can deduce  (c) 
 by applying  the result \cite[Lemma~3.27]{MT1} which is quoted above.
 Finally note that  the only polytopes with $\le 4$ vertices
 are  the simplices $\De_k, k\le 3,$ and  the trapezoid.  Since these all fail condition (c), (d) must hold. 
 Thus these claims hold when $\ell=1$.
   If $\ell > 1$, apply the inductive hypothesis to $F_1$ and use 
   \cite[Prop.~3.17]{MT1}.  

Now consider the face $f$ as a polytope in its own right, and pick 
any two distinct vertices $v_1, v_2$  of $f$.   
Because $\ka$ is generic, every self-equivalence  of  $f$ (resp. $\De$) acts trivially on homology and so belongs to $\Aff_0(f)$ (resp. $\Aff_0(\De)$).
Hence condition (c) implies that no two vertices of $f$  are 
equivalent as vertices of $f$.
It follows easily that they cannot be equivalent in $\De$.  For because
the  self-equivalences $\phi$ of $\De$ belong to $\Aff_0(\De)$ they
are products of the commuting reflections $\rho_i$, where $\rho_i$ interchanges     
 the  pair $F_i, F_i'$ and acts as the identity on all other facets.
 If $v_1,v_2$ are equivalent in $\De$, we may choose $\phi\in \Aff_0(\De)$ which is a product of a minimal number of the $\rho_i$ so that 
 $\phi(v_1)=v_2$.  If  $\phi$  interchanges
$F_i, F_i'$ then, as above,  both $v_1$ and $v_2$ lie in $F_i\cap F_i'$. 
But then $\rho_i\circ \phi$ is a shorter product  that takes $v_1$ to $v_2$, a contradiction.
Thus $\De$ has at least $5$ inequivalent vertices. 
  \end{proof}
  
% The next result shows that the one point blow up of the
%toric manifold $(M,\om)$ has at least two toric structures,
%provided that $\om$ is generic and $M$ is not
%a product of projective spaces    with the product toric structure.

\MS

\NI {\bf Proof of Proposition~\ref{prop:torbl}.}
We must show that the one point blow up of the
toric manifold $(M,\om)$ has at least two toric structures,
provided that $\om$ is generic and $M$ is not
a product of projective spaces    with the product toric structure.
By Corollary \ref{cor:torbl}, it
 suffices to show that if all vertices of $\De(\ka)$ are equivalent for some
generic $\ka$ then $\De$ is a product of simplices.
This will follow from Lemma \ref{le:equivf}  if we show that each facet of $\De$ is equivalent to at least one other facet.
  But given a facet $F$ and a vertex $v'\notin F$, by hypothesis $v'$ is equivalent to every vertex $v\in F$.  Hence  there is $\phi\in \Aff_0(\De)$
that takes $v$ to $v'$.  Therefore $\phi(F)$, which contains $v'$, cannot equal $F$.  
Therefore $F$ is equivalent to at least one other facet, namely $\phi(F)$. 
\QED

  %%%%%%%%%%%%%%%%%%%%%%%%%%%%
  \section{Full  mass linearity}\labell{s:FML}
 %%%%%%%%%%%%%%%%%%%%%%%%%%%%
We begin by improving some results from \cite{MT1} and
then introduce the idea of full mass linearity.   

%%%%%%%%%%%%%%%%%%%%%%%%%%%%%%%%%%%%%%%%%%%%%%%%%%%%%%%%
  \subsection{Some properties of mass linear functions}
  %%%%%%%%%%%%%%%%%%%%%%%%%%%%%%%%%%%%%%%%%%%%%%%%%%%%%%%%

 Let $M$ be a toric $2n$-dimensional 
 manifold with moment polytope $\De$, where
 $$
 \De = \{\xi\in \ft^*: \langle \eta_i,\xi\rangle \le \ka_i, i=1,\dots,N\}.
 $$
  Consider the volume $V(\ka)$ of $\De$ as a function of its support numbers $\ka_i, i=1,\dots,N$.
  The results of Timorin \cite{Tim} 
 show  that the algebra  $H^*(M;\R)$ is isomorphic to
  $\R[\p_1,\dots,\p_{N}]/I(V)$ where 
we interpret $\p_i$ as the differential operator $\frac {\p}{\p \ka_i}$ and
$I(V)$ 
consists of all differential operators with constant coefficients
 that annihilate the polynomial
 $V$; cf. the discussion at the beginning of \S2.6 in \cite{Tim}.
 His argument is the following.   He observes that the
  translational invariance of $V$ implies that
 $\sum \langle\xi,\eta_j\rangle\, \p_jV = 0$ for all $\xi\in \ft^*$.
 Further, he shows that $\p_IV=:V_I$ is the volume 
$$
{\rm vol\,} F_I: = \tfrac 1{k!} \int_{PD(F_I)} \om^k
$$
of the K\"ahler submanifold $\Phi^{-1}(F_I)$ of the face $F_I$.   
(Here $M$ is equipped with the natural 
symplectic form $\om = \om(\ka)$
whose integral over the $2$-sphere corresponding to each edge is the affine length of that edge.)
It follows that
 $\p_I V=0$ whenever $F_I=\emptyset$.   He then shows in 
 \cite[Theorem~2.6.2]{Tim}   that
 these relations generate $I(V)$.   It follows immediately that
  his algebra is isomorphic  to the Stanley--Reisner presentation for $H^*(M)$ described in
  equation \eqref{eq:SR}.

Note that this isomorphism takes $\p_i$ to the 
Poincar\'e dual of the facet $F_i$.  Hence the first Chern class  of $M$ is represented by the operator $\sum_i\p_i$.  More generally, the $k$th
Chern class   is represented by the operator $\sum_{|I|=k}\p_I$.

Consider an element $H\in \ft$.  By taking  inner products, we get an induced function, also denoted $H$, from $\De(\ka)\to \R$.  This is said to be {\it mass linear} if  $H(B_n)$ is a linear function of the $\ka_i$, where $B_n$ is the barycenter of $\De(\ka)$.  Thus there are constants $\ga_i\in \R$ such that $H(B_n) = 
\sum \ga_i\ka_i$.  It is proved in \cite[Lemma~3.19]{MT1} that in this situation
the vector $H\in \ft$  is precisely $\sum \ga_i\eta_i$.
 Thus, if $H$ is mass linear, there are constants $\ga_i$ such that
\begin{equation}\labell{eq:ML}
H(B_n) = 
\sum \ga_i\ka_i,\quad \mbox{and } \;\; H = \sum \ga_i\eta_i.
\end{equation}

If $\mu$ denotes the moment $\int_\De H\,d\Vol$ of $H$ we have
$\mu = H(B_n) V$.    Generalizing Timorin's ideas, we proved the following result in \cite[Proposition~2.2]{MT1}. 

\begin{lemma}\labell{le:tim}
For any $H\in \ft$ the face $F_I$
has volume $V_I: = \p_I V$ and $H$-moment  $\mu_I= \p_I \mu$.
\end{lemma}

Therefore, in the mass linear case we have
$$
\mu_I = H(B_n) V_I + \sum \ga_i V_{{I\less i}}.
$$
The following combinatorial result improves some of the conclusions
of \cite{MT1}.    
We denote by $\vareps_{vj}$  the directed edge 
that starts at the  vertex $v$ and ends transversely to $F_j$.

 \begin{prop}  \labell{prop:FML0} 
 Suppose that  $H\in \ft$ is mass linear and $ H(B_n)=\sum \ga_j \ka_j$.
 Then:
 \begin{itemize}\item[(i)]
$\sum_j \ga_j = 0$.
\item[(ii)] $\sum \ell(\vareps_{vj})\ga_j = 0$ where the sum is over all   directed edges $\vareps_{vj}$, and $\ell(\vareps)$ denotes the affine length of $\vareps$.
\end{itemize}
\end{prop}
\begin{proof}  (i)  Fix a vertex $v: = F_I$.  For each $i\in I$ there is a unique edge $_i\vareps_v = F_{I\less i}$ that starts at $v$ transversely to $F_i$.    Its other endpoint is transverse to a 
unique facet $F_j$ where $j\notin I$.  (Thus $_i\vareps_v = \vareps_{vj}$ in the previous notation.)
For each $j\notin I$  define 
$I(j)\subset I$  to be the (possibly empty) set of $i$ such that the second endpoint of $_i\vareps$  is transverse to $F_j$.
Then the sets $I(j), j\notin I$, form a partition of $I$.  Correspondingly, the sets $J^v(j): = I(j)\cup \{j\}, j\notin I,$ form a partition of $\{1,\dots,N\}$.
Therefore, (i) will follow if we show that
$$
\sum_{i\in J^v(j)}\ga_i = 0, \quad\mbox{ for all } j\notin I.
$$

But this holds by the following calculation.  Fix $j\notin I$ and let
$K(j): = I\cup \{j\}$. 
We first claim that for each $k\in K$, the vertex $F_{K(j)\less k}$ is nonempty exactly if $k\in J^v(j)$.   This is clear if $k=j$.  Otherwise $k\in I$ and $F_{K\less k}= F_{I\less k}\cap F_j$ is nonempty exactly if the second endpoint of 
$_k\vareps_v$  lies on $F_j$, in other words exactly if    $k\in I(j)$.
Using  Lemma \ref{le:tim} and the fact that the intersection of every set  of $n+1$ facets is empty,  we now find that
$$
0=\mu_{K(j)} =  \p_{K}(\sum \ga_i\ka_i\,V) =  \sum_{i\in K} \ga_i \,V_{K\less i} = \sum_{i\in J^v(j)} \ga_i.
$$

Now consider (ii).    Given $K\subset \{1,\dots,N\}$  with $|K|=n+1$, define
$\Ee(K)$ to be the set of  all edges
$F_L$, where $L := L_{s,t}: =  K\less \{s,t\}$ is an edge with endpoints 
 $w_s: = F_{K\less \{s\}}$ and $w_t: = F_{K\less \{t\}}$.  These sets partition
 the set of all edges of $\De$, since for any edge $\vareps$ 
 the set $K(\vareps): = \{i: F_i\cap \vareps\ne \emptyset\}$ has precisely
$n+1$ elements, and $\vareps\in \Ee(K(\vareps))$. 

If  $\Ee(K) \ne \emptyset$, pick any directed edge $\vareps_{vj}\in \Ee(K)$.
Then $K= I\cup \{j\}$, in the language of (i).  Consider any edge $L_{s,t}\in \Ee(K)$.
If $s=j$ then the edge has endpoints $v=w_j$ and $w_t\in F_j$, and so is 
the edge previously called $_s\vareps_v$.  Otherwise $w_s,w_t\in F_j$.  
Observe that $w_s$ and $w_t$
are joined by the edge $F_{(I\less\{s,t\})\cup \{j\}}\in \Ee(K)$. It follows that
the edges $L_{s,t}$ in $\Ee(K)$ form a complete graph.   
Moreover these are the edges  of the dimension $m$ face 
$f: = \cap _{s\in I\less I(j)} F_s$, where $m=|I(j)|$.  Hence this face is a 
simplex, so that all its edges have the same length $\la$.

We need to calculate the sum of $\ell(\vareps_{vj})\ga_j$ over
directed edges.  But this equals the sum of $\ell(\vareps)(\ga_s+\ga_t)$
 over unoriented edges, where $\vareps$ joins $F_s$ to $F_t$. 
We proved in (i) that
$$
\sum_{L_{s,t}\in \Ee(K(j))} (\ga_s+\ga_t) = 0.
$$
Since the sets $\Ee(K(j))$  partition the edges of $\De$, this proves (ii).
\end{proof}

 \begin{rmk}\labell{rmk:comb}\rm  The above proof 
 shows that the coefficients $\ga_i$ of a mass linear function satisfy many enumerative identities, that is, identities that depend only on the combinatorics of $\De$.  Thus the existence of a mass linear function 
 imposes many restrictions on the combinatorics of $\De$.
 For example, if  there is no edge from the vertex $v$ to  $F_j$, then
 the equivalence class $J(j) = I(j)\cup\{j\}$ in (i)  consists only of $\{j\}$, 
 and we conclude  that $\ga_j=0$.  This reproves the result
 in  \cite[Prop.~A.2]{MT1} that every asymmetric facet is powerful.
   Using this, one can immediately deduce that many polytopes have no nonzero mass linear functions $H$.  For example, no polygon with  more than 
 four edges has such an $H$.  
 
 As another example, suppose that $\De$ is any polytope other than a simplex and  blow it up  at one of its vertices $v_0$ to obtain $\De'$.
    Then, because $\De$ has $> n+1$ vertices,
     there is a vertex $w$ in $\De$ that is not connected to $v_0$ by an edge and so is not connected to the exceptional divisor $F_0'$ of $\De'$ by an edge.  
     Therefore if $H$ is a nonzero mass linear function on $\De'$, the coefficient
     $\ga_0$ in the expression $H(B_n(\De'))$ must vanish; in other words the exceptional divisor 
     $F_0'$ is symmetric.
     A similar argument shows that every $H$-asymmetric facet $F_i'$ of $\De'$  (i.e. one with $\ga_i\ne 0$) must meet $F_0'$.    For otherwise, the corresponding facet $F_i$ of $\De$ does not meet $v_0$. Since $\De$ is not a simplex there is an edge $\vareps$ from $v_0$ which does not meet $F_i$.  In the blow up, this edge meets $F_0'$ in a vertex $v'$ which is not joined to $F_i'$ by an edge.  (There is only one edge from $v'$ that does not lie in $F_0'$, namely the blow up of $\vareps$.)   Hence $F_i'$ is not powerful, contradicting our previous results.
 
 This discussion is taken much further in the papers \cite{MT1,MT2}, that classify all mass linear functions on polytopes of dimensions $\le 4$.
 However, rather than focussing on  combinatorial identities these papers
analyze the properties of the symmetric and asymmetric facets.
% In particular, one can interpret the complete graphs found in (ii) in terms of 
% the structure of the  symmetric and asymmetric facets  of $\De$.
 \end{rmk}  
 
  %%%%%%%%%%%%%%%%%%%%%%%%%%%%
  \subsection{Full mass linearity}
 %%%%%%%%%%%%%%%%%%%%%%%%%%%%
     
The following condition was suggested by the work of Shelukhin which is discussed further in \S\ref{ss:shel} below.

\begin{defn}\labell{def:FML}  Let $H\in \ft$.   For each $s=0,1,\dots,n,$ let $$
V^s := \sum_{|I|=n-s}V_I
$$ 
be the sum of the volumes of the faces of dimension $s$, and let $$\mu^s:  
= \sum_{|I|=n-s}\mu_I
$$
 be the sum of the corresponding $H$-moments.
Define $B_s$ to be the center of mass of the facets $F^{s}: = \cup_{|I|=n-s} F_I$.
Thus $B_n$ is the usual center of mass and
 $B_0$ is the average of the vertices.
Then we say that $H$ is {\bf fully mass linear} if 
$H(B_s) = H(B_n)$ for all $s=0,\dots, n-1$.
\end{defn}

Note the following points.
\MS

\NI $\bullet$
 Since $B_0$ 
is clearly a linear function of the support numbers 
$\ka_i$, every fully mass linear function is mass linear.
\MS

\NI $\bullet$
 Every inessential function is fully 
mass linear since the barycenters
 $B_s$
must lie on all planes of symmetry of $\De$, i.e.
 they are invariant under the action of elements in $\Aff_0(\De)$. 
\MS

\NI $\bullet$  We explain in \S\ref{ss:shel}   Shelukhin's argument
that the quantities $H(B_s) - H(B_n)$ are values of certain real-valued  characteristic
 classes for Hamiltonian  bundles  with fiber $(M_\De,\om_\ka)$.
It follows that   {\it the function $H$ is fully mass linear whenever
$\La_H$ has finite order in} $\pi_1\bigl(\Ham(M_\De,\om_\ka)\bigr)$.
   In fact, there is one such characteristic class $I_\be$ for each product 
   $c_\be$ of Chern classes on $M$.  
   However, 
 as we show in Corollary \ref{cor:cohom}, the vanishing 
 of these classes  $I_\be$ gives no new information.
%However, the full mass linearity condition is easier to understand in terms of 
%the combinatorics of $\De$.   
\MS

We continue our discussion 
by explaining precisely what full mass linearity means.

\begin{lemma}\labell{le:FML} Let    $H\in \ft$ be mass linear with $H(B_n) = \sum\ga_i\ka_i$.  Then $H(B_{n-r}) = H(B_n)$ 
exactly if the identity
$$
(*_{r})\qquad\qquad \sum_{i,J: i\in J, |J|=r-1} \ga_i V_{J} = 0,
 $$
holds, where we interpret $(*_1)$ to be the identity $\sum\ga_i=0$.
  In particular, $H$ is fully mass linear exactly if $(*_r)$ holds for 
$r=1,\dots,n$.
\end{lemma}
\begin{proof} 
 First consider   $\mu^{n-1}$.
By Timorin \cite{Tim} we have $\mu^{n-1}  = \sum_i \p_i\mu$.
Therefore, because $H$ is mass linear,
\begin{eqnarray*}
H(B_{n-1}) V^{n-1} = \mu^{n-1}  &=& \sum_i \p_i\mu\\
&=& \sum_i H(B_n) V_{i}+ \sum_i \p_i\bigl(\sum \ga_j\ka_j\bigr)V\\
&=& H(B_n) V^{n-1} + (\sum \ga_i) V.
\end{eqnarray*}
This proves the case $r=1$ of the first statement.
Note also that because 
 $\sum \ga_i = 0$ by Proposition \ref{prop:FML0} (i), we always have $H(B_{n-1}) = H(B_n) $.
 
More generally, since $\sum \ga_i = 0$, we have
\begin{eqnarray*}
 H(B_{n-r}) V^{n-r} =\mu^{n-r} &=& \sum_{|I|=r} \p_I\bigl(\sum \ga_j \ka_j V\bigr)\\
&=& H(B_n) V^{n-r} + \sum_{i\notin J, |J|=r-1} \ga_i \p_J V\\
&=&  H(B_n) V^{n-r} +  (\sum_i \ga_i)\sum_{|J|=r-1}\p_J V - \sum_{i\in J, |J|=r-1}\ga_i \p_J V\\
& = &  H(B_n) V^{n-r} -\sum_{i\in J, |J|=r-1} \ga_i V_{J}
\end{eqnarray*}
Therefore  we see that
 the identity  $H(B_{n-r}) = H(B_n)$ holds 
exactly if
$(*_r)$ holds.  This proves the first statement.
The second is clear.
\end{proof}

\begin{rmk}  \rm  The identity $(*_{n+1})$
is  $\sum_{i\in J, |J|=n} \ga_i V_{J} = 0$, which is equivalent to saying that $\sum_i N_i\ga_i = 0$ where  $N_i$ is
the number  of vertices in the facet $F_i$.  But this holds for 
all mass linear functions, as one can see by computing
$0=\sum_{|I|=n+1} \p_I\bigl(\sum \ga_j \ka_j V\bigr)$ as above.
Another way to calculate this is to think of it as a sum over directed
edges, namely 
$\sum_{\vareps_{vj}} \ga_j$.  We proved that this sum vanishes
 in the course of proving
part (ii) of  Proposition \ref{prop:FML0} since the lengths turned out to be irrelevant.
\end{rmk}

\begin{prop}\labell{prop:FML}  
The following conditions are equivalent:
\begin{itemize}  \item 
 $H$  is mass linear;
 %and $H(B_n) = \sum \ga_i\ka_i$ where  $\sum \ga_i = 0$;
\item $H(B_n) =H(B_0)$;
\item  $H(B_n) =H(B_{n-1}) = H(B_0)$.
\end{itemize}
\end{prop}
\begin{proof} 
 If $H$ is mass linear, then we saw in 
the proof of Proposition \ref{prop:FML} that  
$H(B_n) = H(B_{n-1})$  because
 $\sum \ga_i = 0$.   
 Further,  the difference between $
 \mu^0 = H(B_0)V^0$ and  $H(B_n)V^0$ is  
 $$
 \sum_{i\in J,|J|=n-1}\ga_iV_J = \sum_{i\in J,|J|=n-1}\ga_j \ell(F_J)
 $$
 But we saw in Proposition \ref{prop:FML0}(ii) that this sum vanishes.
 Hence the first condition implies the second and third.
 
 But we noted earlier that $H(B_0)$ is a linear function of the $\ka_j$.  Hence the second condition implies the first.  
\end{proof}

\begin{rmk}\labell{rmk:FML}\rm  (i)  This argument shows that the identity $H(B_n) = H(B_{0})$ implies $H(B_n) = H(B_{n-1})$.    Thus if $H$ is mass linear these three points always  lie on the same level set of $H$.  In contrast, 
Shelukhin \cite{Shel} showed in the monotone case that   the three points $B_n, B_{n-1}$ and $B_0$ are collinear.
\MS

\NI (ii) The $r$th equation in Lemma \ref{le:FML} corresponds to a condition on $\mu^{s}, $ where $s=n-r$ that we calculate assuming that $H$ is mass linear.   Therefore this equation 
 is not equivalent to the fact  that $H(B_{n-r}) = H(B_n)$.
 In fact, we give an example in  Remark \ref{rmk:geom} (ii) below 
 showing  that the identities $(*_r), r=1,\dots,n$ do not by themselves imply mass linearity.   
 \end{rmk}

\begin{cor}\labell{cor:FML} Suppose that the mass linear function $H$ has coefficients $\ga_i$ as in Equation \eqref{eq:ML}.  Then $H$  is fully mass linear  exactly if  $\sum \ga_i \p_i^k V = 0$ for all $k=1,\dots,n$.
\end{cor}
\begin{proof}   These identities are equivalent to
$(*_r)$, $r=1,\dots,n$ because the functions $\sum x_i^k$ form a basis for the symmetric polynomials over $\Q$.
\end{proof}

\begin{rmk}[Geometric interpretation of equations $(*_{2})$.]\labell{rmk:geom}
\rm  (i)  Inessential mass linear functions $H$ are generated by
vectors $\xi_H\in \ft^*$ with the property that the facets
$\{F_i: \langle\eta_i,\xi_H\rangle\ne 0\}$  are all equivalent.  
In this case, we saw that  $H(B_n(\ka)) = 
 \sum  \langle\eta_i,\xi_H\rangle \ka_i$, as well as
$H = \sum  \langle\eta_i,\xi_H\rangle \eta_i$.
Moreover, if $H$ is elementary, i.e. of the form $\eta_i-\eta_j$,
 then by \cite[Lemma~3.4]{MT1} there is an affine reflection symmetry of $\De_H$ that interchanges the two facets   
$F_i$ and $F_j$  preserving the transverse vector $\xi_H$.
 
We claim that a very similar statement holds for mass linear functions $H$
that satisfy $(*_2)$.  In other words, for each such function there is a vector $\xi_H\in \ft^*$
such that
\begin{equation}\labell{eq:ML2}
\ga_i =  \langle\eta_i,\xi_H\rangle \;\;\mbox{ where } \;
H(B_n(\ka)) = 
 \sum  \ga_i \ka_i.
\end{equation}
To see this, observe that
 equation  $(*_2)$: $\sum_i\ga_i V_i = 0$  says that the operator $\sum \ga_i \p_i$ is in the annihilator $I(V)$.  Timorin showed that $I(V)$ is generated by additive relations of the form
 $
\sum  \langle\eta_i,\xi\rangle \p_i=0$  where $\xi\in \ft^*$,
as well as some multiplicative relations $\p_I = 0$. 
Since $\sum \ga_i \p_i$ is linear, it has to correspond to some vector 
$\xi_H\in \ft^*$.  Note that the first part of Equation \eqref{eq:ML2}
shows that $\xi_H$ must be parallel to all symmetric facets. However, it is not clear whether there is further geometric significance to this vector.

This observation explains the condition $\sum \ga_i\, a_i = 0$ in
Lemma \ref{le:prod} below.  For in this case $\xi_H = -(\ga_1,\dots,\ga_k,0) \in \R^{k+1}\equiv \ft^*$ while the two  facets with conormals
$\eta_{n+1} = (0,\dots,0,1)$ and $\eta_{n+2} = (-a_1,\dots,-a_k,1)$
are symmetric.\MS

\NI (ii)  In general, one cannot reduce the mass linearity condition for $H: = \sum_i\langle \eta_i,\xi_H\rangle \eta_i$  to any obvious condition on 
$\xi_H$.  Consider for example the $\De_1\times \De_1$-bundle over $\De_1$ with conormals 
$$
\eta_1 =  -e_1, \eta_2= -e_2, \eta_3 =e_1,\eta_4=e_2, \eta_5: = -e_3, \eta_6: = e_3-v,
$$
 where $v=(a_1,a_2,0)$ as in Lemma \ref{le:prod}.  For generic $(a_1,a_2)$ (i.e. $a_1a_2\ne 0$, and $a_1-a_2\ne 0$),  this has just one pair of equivalent facets, namely the base facets $F_5, F_6$.  Since the other facets are neither pervasive nor flat, \cite[Theorem~1.10]{MT1} implies that
$\De$ has no mass linear functions for which $F_5, F_6$  are symmetric.  On the other hand, if $\xi_H: = (-a_2,a_1,0)$  we get
$H = a_2(\eta_1 -\eta_3) -a_1( \eta_2 - \eta_4)$. 
So this $H$ satisfies
$(*_1)$, and it satisfies $(*_2)$ by construction.
One can easily check that  $(*_3)$ holds.  
  Thus, by Proposition \ref{prop:FML0}, $H$ satisfies all the identities in Lemma \ref{le:FML}, but  it is not mass linear.
%\MS
%
%
%\NI  (iii) It seems possible that the conditions for 
%$H$ to be fully mass linear are stronger than those for $H$ to 
%be mass linear.  However, no examples are known.  
\end{rmk}

  %%%%%%%%%%%%%%%%%%%%%%%%%%%%
  \subsection{Examples}\labell{ss:ex}
 %%%%%%%%%%%%%%%%%%%%%%%%%%%%

  We now describe one of the basic examples from \cite{MT1,MT2}.
 Suppose that $\De\subset \R^k\times \R$ is a $\De^k$-bundle over $\De^1$  
% as in Definition \ref{def:bun}. Thus its shape is determined
% by the  constants $a_1,\dots,a_k$.  It will be convenient  
% to set $a_{k+1}: = 0$.
with conormals
\begin{gather}\labell{eq:prod}
\eta_i = -e_i, i=1,\dots, k,\quad \eta_{k+1} = \sum_{i=1}^k
e_i ,\\\notag \eta_{k+2} = -e_{k+1}, \quad \eta_{k+3} = e_{k+1} + \sum_{i+1}^k a_i e_i;
\end{gather}
cf. Definition \ref{def:bun}.  Thus $\De$ is determined 
by the vector $A: =  (a_1,\dots,a_k)$. %= \sum_{i=1}^{k+1} \ga_i\eta_i$,
For  convenience  
we later set $a_{k+1}: = 0$.

\begin{lemma}\labell{le:prod}  With $\De$ as in equation \eqref{eq:prod},
 the function
$H = \sum_{i=1}^{k+1} \ga_i\eta_i$ is fully mass linear  exactly if it is mass linear, which happens exactly if  
$$
\sum \ga_i=0, \;\;\;\mbox{ and } \;\;\; \sum \ga_i\,a_i = 0.
$$
\end{lemma}
\begin{proof} It is easy to check that
the volume function of $\De$  is  
$$
V(\ka) = \tfrac1{k!} h\,\la^k - \tfrac1{(k+1)!}\bigl(\sum a_i\bigr)\,\la^{k+1},
$$
where 
$$
h = \ka_{k+2} + \ka_{k+3} + \sum_{i\le k+1} a_i\ka_i, \;\;\;
\la = \sum_{i=1}^{k+1}\ka_i.
$$
%Details may be found in \cite[\S4]{MT2}. Moreover, we
Moreover, one can show
by direct calculation  that
$H = \sum_{i=1}^{k+1} \ga_i\eta_i$ is  mass linear on $\De$ exactly if 
$\sum \ga_i=0$ and $ \sum \ga_i\,a_i = 0$.   The case $k=3$ is worked out in detail in
\cite[Proposition~4.6]{MT1}.  The general case is similar; details will appear in 
\cite[\S4]{MT2}.

Therefore we need to show that these two conditions imply that
$H$ is fully mass linear.  By Corollary \ref{cor:FML} it suffices to see that $\sum\ga_i\p_i^m V= 0$ for all $m$.  
This is an easy calculation.  Note also that when $m=1$ this condition says that
$\sum \ga_ia_i=0$ and is equivalent to the statement that $H(B_{n-2}) = H(B_n)$. \end{proof}  

\begin{cor}  Every mass linear function on a polytope of dimension $d\le 3$ is fully mass linear.
\end{cor}
\begin{proof}  This is an immediate consequence of Proposition \ref{prop:FML} when 
$d=2$, and is anyway clear because all mass linear functions in $2$ dimensions are inessential, and hence fully mass linear.  We showed in \cite{MT1} that when $d=3$
the only essential mass linear $H$ occur on polytopes that are $\De_2$ bundles over $\De_1$, and (modulo adding an inessential function) are of the form considered in Lemma \ref{le:prod}. 
  Hence the result follows from Lemma \ref{le:prod}.
\end{proof}

\begin{rmk}\labell{rmk:4}\rm
 In dimension $4$, it is easy to check that a mass linear function $H$
%\begin{lemma}\labell{le:dim4}    Let $H$ be a mass linear function on a $4$-dimensional
%polytope $\De$ with $H(B_n) = \sum\ga_i\ka_i$.  
%Then 
 is fully mass linear if:\begin{itemize}
\item
 it is geometrically generated; i.e. there is  vector $\xi_H\in \ft^*$ such that $\ga_i = \langle \eta_i,\xi_H\rangle$ for all $i$;  and

\item $\sum_{i\ne j, i,j \in \Aa} \ga_i V_{ij} = 0$ where $\Aa=\{i: F_i \mbox { is asymmetric}\} = \{i:\ga_i\ne 0\}$.
\end{itemize}
The pairs $(\De,H)$ where $\De$ has dimension $4$ and $H$ is essential
are classified in \cite{MT2}.  It appears that in all cases $H$ is fully mass linear. 
It would be interesting to find a more conceptual proof; the classification in \cite{MT2} is too complicated to transfer easily to higher dimensions.
\end{rmk}
  \subsection{Mass linearity and characteristic classes}\labell{ss:shel}
 %%%%%%%%%%%%%%%%%%%%%%%%%%%%

 We now explain Shelukhin's approach to mass linearity.
 Every Hamiltonian bundle  $P\to S^2$ with fiber $(M,\om)$  
  carries
a canonical extension $u\in H^2(P;\R)$ of the class of the symplectic form on $M$ called the {\it coupling class}.\footnote
{This is the unique extension such that $\int_{P} u^{n+1} = 0$.}
One also considers 
the vertical Chern classes $c^{\Vertt}_{n-s}\in H^*(P)$,
which are just the ordinary Chern classes 
of the tangent bundle to the fibers of $P\to S^2$.
Using this data one can define a homomorphism
$\pi_1\bigl(\Ham (M,\om)\bigr) \to \R$ by integrating a product of some vertical Chern classes with a suitable power of $u$ over $P$.   For example, we define
$I_s$  by integrating 
$c_{n-s}^{\Vertt} u^{s+1}$. \footnote
{
These characteristic classes were first
 defined in  \cite{LMP}; see also \cite[\S3]{KM}.}

If the element $\La_H\in \pi_1\bigl(\Ham (M,\om)\bigr) $ is toric, then as we saw above 
$M_H$ is toric.  Moreover, for each $s= 0,\dots,n-1$,
the class  $c^{\Vertt}_{n-s}$
 is Poincar\'e dual to
$F_H^{s}$, the union of the faces of $\De_H$ of dimension $s+1$ and transverse to the fiber, i.e the union of the prolongations to $\De_H$ of  all faces of $\De$ of  dimension $s$.

Shelukhin showed in \cite[Thm.~4]{Shel}
 that $H$ is fully mass linear if and only if the corresponding loop 
$\La_H\in \pi_1\bigl(\Ham (M,\om)\bigr) $ is in the kernel of the 
homomorphisms $I_s$, for $0\le s< n$. In fact, by finding a nice representative for  the coupling class $u$ in terms of the normalized Hamiltonian $H - H(B_n)$,
 he showed that
\begin{equation}\labell{eq:I0}
I_{s}(\La_H) = {\it const}\, \int_{F^s} \bigl(H-H(B_n)\Bigr) d\Vol
= {\it const}\, \bigr( H(B_s) - H(B_n)\bigr) V^s.
\end{equation}
This motivated Definition \ref{def:FML}:  since our work on mass linear functions is  primarily aimed at  understanding
 the kernel of the map $\pi_1(T)\to \pi_1\bigl(\Ham (M,\om)\bigr) $, fully mass linear functions are really more relevant to us than mass linear ones.  However, the examples in the previous section show that mass linearity 
 seems to be the most crucial part of the full condition.
  
 More generally, given any tuple $\be: = (\be_1,\dots,\be_n)$  with $|\be|: = \sum i \be_i \le n+1$, set $c^{\Vertt}_\be :=
 \prod (c^{\Vertt}_i)^{\be_i}$ and define
 \begin{equation}\labell{eq:be0}
 I_\be (H) = \int_{M_H} c^{\Vertt}_\be u^{n+1-\be}.
 \end{equation}
 Shelukhin also observed that $ I_\be (H)$ must vanish if
 $\La_H$ has finite order in $\pi_1(\Ham)$.  If  
  $c_\be: = \prod (c_i)^{\be_i}$  is represented by the weighted sum $
 \sum_{|I|=|\be|} m_I F_I$  of faces of $\De$, then as above
 \begin{equation}\labell{eq:be}
  I_\be (H) = {\it const} \sum m_I\bigl(H(B_{F_I}) - H(B_n)\bigr),
\end{equation}
  where $B_{F_I}$ is the barycenter of $F_I$.  
  
  \begin{lemma} \labell{le:tot} If $H$ is fully mass linear then
   $I_\be(H) = 0$ for all $\be$.
  \end{lemma}
 
  We prove this in Corollary \ref {cor:cohom}; it is a consequence of our cohomological description of mass linearity.
  %We  call functions $H$ for which
  %all $I_\be(H)$ vanish {\it totally mass linear}.
 
 Some of these classes always vanish by the standard ABBV localization formula.
A particularly easy case is when
$\be=c_1c_{n}$. Then 
 $$
 I_\be(H) = \int_{F_H^{n}} c_1^{\Vertt}
 $$
 is the integral of $c_1^{\Vertt}$ over the edges of 
$\De_H$ that do not lie in any fiber.  Modulo a constant, 
this is simply the sum of the isotropy 
 weights of $H$  at the vertices of $\De$ and so always vanishes.
\footnote
{
Here is a brief proof:  
Because $I_{\be}(H)$ is linear in $H$ it is enough to prove this
for a set of integral $H$ whose rational span includes $\ft_\Z$. 
Therefore we can assume that the critical points of $H$ are just 
the vertices of $\De$, and
that at each vertex the weights are pairwise linearly independent.
Then the set of points in $M$ with nontrivial stablizer is a union 
of $2$-spheres; each has exactly two fixed points with opposite weights.
 See Pelayo--Tolman \cite[Thm.~2,Lemma~13]{PT} for a much more 
precise version of this result that uses the ABBV localization in its proof.}
(As explained by Shelukhin, these are special cases of some 
vanishing results for Futaki invariants.)

\begin{rmk}\labell{rmk:TML}\rm  Formula \eqref{eq:be} holds for all ways of representing 
the class $c_\be$ as a sum of facets.   This gives yet more identities 
that have to be satisfied by fully mass linear functions. But many of these will be automatically satisfied.  For example,  if two facets $F_1$ and $F_2$ are homologous, then there is an affine self-map of $\De$ that interchanges them (cf. the discussion after Definition \ref{def:equiv}).  Hence $H(B_{F_1})$ is a linear function of $\ka$ if and only if 
$H(B_{F_1\cup F_2 })$ is.    Similarly if two faces $F_I$ and $F_{I'}$ are homologous there may well be  an affine self-map   that interchanges them.  However, in the absence of such we might get new information.  This could be combined with an analysis of the asymmetric and symmetric facets
considered   in \cite{MT1,MT2}.
\end{rmk}

%\begin{example}\labell{ex:TML}
%\rm   Let  $M$ be the one point blow up of $\C P^3$ with polytope $\De$   given by $$
%x_1,x_2,x_3\ge 0, \;\;\sum x_i\le\la, \;\;x_3\le h.
%$$
%Let $a\in H^2(M)$ be dual to the facet $x_3=0$ and $e$ be dual to $x_3=h$, so that $a^3 = pt = -e^3$. Because $c_1$ is represented by the sum of the facets (cf. \cite{DJ}), we have $c_1= a + 3(a-e) + e = 4a-2e$.  Therefore $c_1^2 = 16 a^2 +4e^2$.  On the other hand, $c_2$ is represented by the sum of the intersections of pairs of facets, so that  $c_2 = 3 a^2 + 3(a^2-e_3) + 3e^2 = 6a^2.$
%Therefore, if $H$ is fully mass linear, 
% the value of $H$ on the center points $B_\vareps$ of the edges representing
%$a^2$ and $e^2$, and hence also those representing $a^2-e^2$, must be linear functions of $h,\la$.   In other words $H(B_\vareps)$ must be a linear 
%function of $\la,h$ for all edges $\vareps$.   This is not an interesting example, because the center points themselves vary linearly with $h,\la$.  But  in higher dimensions it is possible that
%the full mass linearity condition would give more information 
%than mass linearity.  For example in $4$-dimensions  the vanishing of the homomorphisms  $I_{c_1^2}$ and 
%$I_{c_2}$ would give information 
%on the barycenters of sums of $2$-faces.  
%\end{example} 

 %%%%%%%%%%%%%%%%%%%%%%%%%%%%%%%%%%%%%%%%%%%%%%%%%%%%%%%%
\subsection{A cohomological interpretation of mass linearity}\labell{ss:cohom}
 %%%%%%%%%%%%%%%%%%%%%%%%%%%%%%%%%%%%%%%%%%%%%%%%%%%%%%%%

We saw in Lemma~2.5 of \cite{MT1} that the set of mass linear functions
$H\in \ft$
forms a rational subspace of $\ft$, and hence 
 is generated by elements of the integer lattice $\ft_\Z$ of $\ft$.  Hence we will restrict attention here to  $H\in \ft_\Z$.  Each such $H$ exponentiates to a circle subgroup $\La_H$ of the Hamiltonian group of
 the toric manifold $(M_\De,\om)$, and as before, we denote by $M_H$
 the corresponding fibration over $S^2$ with fiber $M$ and clutching map $\La_H$.  
  In this section we describe what it means for $H$ to be mass linear in terms of the cohomology algebra of $M_H$.

We now investigate the volume function $V^H$ of $\De_H$.
Note that $\De_H\subset \ft^* \times \R$ has 
$N$ facets $F_j^H$ corresponding to the $F_j$ in $\De$ with conormals $(\eta_j,0)$,
and two other facets $F_{N+1}, F_{N+2}$ with conormals
 $\eta_{N+1} = (0,-1), \eta_{N+2} = (0,1) + H$; cf.   Thus, because we may write $H = \sum_{i\le N} \ga_i \eta_i$ we have
\begin{equation}\labell{eq:etaN+1}
\eta_{N+1} +\eta_{N+2} -  
\sum_{i\le N} \ga_i \eta_i=0.
\end{equation}
Further the top facet is given by points 
$(\xi,t)\in \ft^*\times \R$ such that
$t+H(\xi) = \ka_{N+2}$.
%Write $V^H: = $ volume of $\De^H$. This is equal to the $H$-moment $\mu_H$.
%Hence, by mass linearity,
The volume $V^H$ of $\De_H$ is therefore
%%
%\begin{equation}\labell{eq:VH}
$$
V^H =\int_{\De} \Bigl(\int_{-\ka_{N+2}}^{\ka_{N+1}-H(\xi)} dt\Bigr) d \Vol(\xi) 
= \Bigl(\ka_{N+1}   + \ka_{N+2} - H(B_n)\Bigr) V,
$$
%\end{equation}
%%
where $B_n$ is the center of gravity of $\De$.

By Timorin, 
$$
H^*(M_H) \cong \R[\p_1,\dots,\p_{N+2}]/I(V^H) =: R_H
$$ 
where 
we interpret $\p_i$ as the differential operator $\frac {\p}{\p \ka_i}$ and
$I(V^H)$ 
consists of all differential operators that annihilate the polynomial
 $V^H$.
  The multiplicative relations  in $I(V^H)$ are 
$\p_{N+1} \p_{N+2} = 0$ 
together with all  multiplicative relations $ \p_I = 0$ for $V$.
Since there is also a new additive relation $\p_{N+1} -\p_{N+2} = 0$,
 we will from now on set  $\ka_{N+2} = \ka_{N+1}$ and use the relation $\p_{N+1}^2 = 0$.  Therefore we take $V^H$ to be
\begin{equation}\labell{eq:VH}
V^H =\Bigl(2\ka_{N+1} - H(B_n)\Bigr) V,
\end{equation}
%%
%Therefore if we define $\p'_j : = \p_j + \ga_j\p_{N+1}$ the $\p_j'$ satisfy 
%the  additive relations for $V^H$.  

\begin{rmk}\labell{rmk:coeff}\rm
In the next theorem we must be careful about the coefficients.
In order for  $\La_H$ to be a circle action,  we assumed that
$H\in \ft$ is integral.  However, in the mass linear case
 this does not mean that 
the coefficients  $\ga_i$ in the expression $H(B_n) = 
\sum \ga_i\ka_i$ are integers.   For example,  if $\De=\De_1\subset \ft^*=\R$ is the $1$-simplex with conormals $\eta_1=-1, \eta_2=1$, and if
$H=\eta_2\in \ft$, then $H(B_1) = -\frac 12 \ka_1 + \frac 12 \ka_2$.
Correspondingly, $\La_H$ is the rotation of  $S^2= M_\De$ by one full turn, with order $2$ in $\pi_1\bigl(\Ham(S^2,\om)\bigr)$.
In fact, we prove in   \cite[Prop.~1.22]{MT1} that the loop $\La_H$ contracts in $\Ham(S^2,\om)$ only if the $\ga_i\in \Z$.
 It follows that if $\La_H$ has finite order $m$ in $\pi_1\bigl(\Ham(M,\om)\bigr)$,
 then the numbers $m\ga_i$ are all integers.
Note also that the $\ga_i$ are always rational because,
as we point out in \cite[Rmk.~2.4]{MT1},  
the polynomial functions $V(\ka)$ and $\mu(\ka)$ have rational coefficients.

 In Theorem \ref{thm:cohom} below, we consider cohomology with coefficients $\R$.  However, the isomorphism $\Psi$ (if it exists) is rational, and it  induces an isomorphism on
 integral homology exactly if the the coefficients $\ga_i$ are integers.
 (Note that $H^*(M;\Z)$ is torsion free  when $M$ is a toric symplectic manifold.)  Note also that $\Phi$ induces the identity map on the cohomology $H^*(M)$ of the fiber.
\end{rmk}

\begin{thm}\labell{thm:cohom}
Let $(M,\om,T)$ be a toric manifold with moment polytope $\De$, and let 
$H\in \ft\less \{0\}$.  Let $M\to M_H\to S^2$ be the corresponding bundle.
\SSS

\NI
 {\rm (i)}  
The function $H$ is mass linear on $\De$ with $H(B_n) = \sum \ga_i\ka_i$
if and only if
there is an algebra isomorphism 
$$
\Psi:
H^*(S^2)\otimes H^*(M) \equiv \Bigl(\R[z]/z^2\Bigr) \otimes \Bigl(\R[\p_1,\dots,\p_N]/I(V)\Bigr)\;\to\;
H^*(M_H) 
$$
that is compatible with the fibration structure on $H^*(M_H)$, i.e. 
if we identify   $H^*(M_H)$ with the algebra $\R[\p_1,\dots,\p_{N+1}]/I(V^H)$
as above then there are constants $\al_i$ such that
$$
\Psi(z) =
\p_{N+1},\;\; \mbox{ and }\;\;
\Psi(\p_i) =\p_i': =  \p_i+\al_i \p_{N+1}\in I(V^H).
$$  
\MS

\NI {\rm (ii)}  If $H$ is mass linear, then
 it is 
 fully mass linear exactly if  $\Psi$
 takes the Chern classes $c_s^M$ in $H^0(S^2)\otimes H^*(M) $ to the 
vertical Chern classes $c_s^{\Vertt}$ in $H^*(M_H)$ for all $s=1,\dots,n$.
\MS

\NI {\rm (iii)}  If $H$ is mass linear then $\Psi(c_1^M) = c_1^{\Vertt}$
and  $\Psi(c_n^M) = c_n^{\Vertt}$.
\end{thm}

\begin{proof}   Suppose first that $H$ is mass linear.
Then, by equation \eqref{eq:ML}, there are constants $\ga_i$ such that $H= \sum \ga_i\eta_i$ and $H(B_n) = \sum \ga_i\ka_i$.
 Since $\De_H$ is combinatorially equivalent to a product, it follows from the Stanley-Reisner presentation for  $H^*(M_H)=:R_H$  that this algebra is additively isomorphic to 
a product.\footnote
{
In fact this is true for 
 all Hamiltonian bundles over $S^2$ by \cite{LMP,Mcq}.}
 By Equation (\ref{eq:etaN+1})
the  additive relations for $V^H$ are 
$$
0 = \sum_{j\le N} \langle \xi_i,\eta_j\rangle \p_j + \langle \xi_i,
\eta_{N+2}\rangle \p_{N+1} 
= \sum_{j\le N} \langle \xi_i,\eta_j\rangle \bigl(\p_j + \ga_j\p_{N+1}\bigr)
$$ 
where $\xi_i$ runs over a basis for $\ft^*$. Therefore, if we
take $\al_i = \ga_i$ for all $i$, the map $\Psi$ defined in (i) is an additive homomorphism.
Therefore it remains to check that 
the relations $\p_I=0$ that generate the 
multiplicative relations in $I(V)$ are taken by $\Psi$ to  relations
$\p'_I$ in $I(V^H)$.

To see this, note that
 $V_I=0$ iff $F_I=\emptyset$, while $F_I=\emptyset$ implies   $\mu_I = 0$.
 Therefore, because $\mu=(\sum \ga_i\ka_i)V$,  for such $I$ we have
 $$
0=\mu_I = H(B_n)V_I +
 \sum_{i\in I}\ga_i V_{I\less i} =\sum_{i\in I}\ga_i V_{I\less i}.
$$
Hence
\begin{eqnarray*}
\prod_{i\in I} \p_i' V^H &=& 
\prod_{i\in I}\bigl(\p_i+ \ga_i\p_{N+1}\bigr)\bigl(2\ka_{N+1}-H(B_n)\bigr) V\\
& =&
\bigl(2\ka_{N+1}-H(B_n)\bigr) V_I + \sum_{i\in I} (2\ga_i - \ga_i) V_{I\less i}\\
&=& \sum_{i\in I} \ga_i V_{I\less i} = 0,
\end{eqnarray*}
as required.
  
Therefore there is an algebra homomorphism $\Psi: H^*(S^2)\otimes H^*(M)\to H^*(M_H)$.
By construction, its composition  with the restriction map $H^*(M_H)\to H^*(M)$ 
is surjective.  Therefore, by 
 the Leray--Hirsch theorem, it is an isomorphism.

%Moreover these relations have the form
%$F_I^H + R_I'=0$ where $F_I^H=0$ is the Stanley 
% because if this is the case then  the existence of  any other linearly independent relation in $H^*(M_H)$
% would make its rank too small.

Conversely, suppose that
 $$
\Phi:  H^*(S^2)\otimes H^*(M)\to H^*(M_H)
$$
is an isomorphism of algebras that is compatible with the fibration, i.e. 
its restriction to the fiber $H^*(M)$ is
the  identity and it takes the generator of $H^2(S^2)$ to the pullback of
this class in $H^2(M_H)$.
  We must show that $H$ is mass linear.

Let us think of   the symplectic class $[\om] = [\om_\ka]$ on $M$ 
as a function of the 
the support numbers $\ka$ of the polytope $\De=\De(\ka)$.  
In terms of the chosen isomorphism $H^2(M)$ with the degree $2$ part of
the algebra $\R[\p_1,\dots,\p_N]/I(V)$, we may write $[\om_\ka] = \sum 
\ell_i(\ka)\p_i$ where the coefficients $\ell_i(\ka)$ are linear functions of
$\ka$.   
Similarly, the symplectic class $[\Om_\ka]$ (which is determined by 
positions of the facets  of the polytope $\De_H(\ka)$) is a linear function of $\ka$.

Because $\Phi$ restricts to the identity on the fiber and
is compatible with the identity map on the base, 
the induced map on $H^{2n+2}$ preserves the integer lattice, and hence
preserves the  cohomological fundamental class.  Therefore there is a well defined
%(Note also that, because the integral cohomology of a toric 
%manifold has no torsion, the integral cohomology embeds in the real cohomology.)
 the Poincar\'e dual isomorphism 
$\Phi_*: H_*(S^2)\otimes H_*(M)\to H_*(M_H)$ that
 takes the (homology) fundamental class
 $[S^2\times M]$ to $[M_H]$. 
 Further, 
%Let $Z \in H_2(M_H)$ be
 the image $Z: = \Phi_*([S^2])$ of the 
fundamental class of $S^2$
is independent of $\ka$.  Hence 
the above remarks imply that
$$
\int_Z \Om_\ka = L(\ka)
$$
is a linear function of $\ka$.

Now observe that
 the volume $V^H$ is a cohomological invariant of $M_H$: up to a constant,
  it is obtained by evaluating
$(\Om_\ka)^{n+1}$  on the fundamental class in
 $H_{2n+2}(M_H)$.  Thus we can evaluate $V^H$ in 
 the product algebra.  But here it is just the product of the area of $[S^2]$ 
 (with respect to $\Phi^{-1}[\Om_\ka]$)
  with the volume $V$ of $M$. 
Since the area of $[S^2]$ is $\int_Z[\Om_\ka]$  it follows that
$V^H$ has the form $L(\ka_i) V$ where $L$ depends 
linearly on the $\ka_i$ as we saw above.
Because, as we noted in Remark \ref{rmk:coeff}, 
the functions
$V^H$ and $V$ have rational coefficients,  
the coefficients of $L$ must also be rational.
But we saw in 
 Equation \eqref{eq:VH} that $V^H = (2\ka_{N+1} - H(B_n)) V$.
 It follows that $H(B_n)$ is a linear function 
of the $\ka_i$ with rational coefficients.  
This completes the proof of (i).

%For any $\ka$
%the coupling class $u_\ka\in H^2(M_H)$ 
%  is defined to be the unique extension of $[\om_\ka]$ such that $u_\ka^{n+1}=0$.
 % Hence $u_\ka = \Phi([\om_\ka])$ where we identify $[\om_\ka]\in H^2(M)$ 
 % with the corresponding class in $H^0(S^2)\times H^2(M)$.
  %Moreover, the symplectic class $[\Om_{\ka,H}]$ on $M_H$ has the form $u_\ka + 
%\la_\ka\pi^*(a)$,
 % where $\pi:M_H\to S^2$ is the projection, $a\in H^2(S^2)$ is the 
%canoncial generator and $\la_\ka$ is a suitable function of $\ka$.
%  Thus $[\Om_{\ka,H}] = \Phi(\la_\ka a + [\om_\ka])$.
  \MS

Now consider (ii).
The Chern classes $c_s$ of $M$ are Poincar\'e dual 
to the classes in $H_{2n-2s}$ represented by the 
face sums $\sum_{|I|=s} F_I=: F^{s}$.
Thus they are represented in the algebra $\R[\p_1,\dots,\p_N]/I(V)$
by the differential operator $\sum_{|I|=s} \p_I$.
These same operators also represent the vertical Chern classes 
of the trivial bundle $S^2\times M$.

Next observe that the  vertical Chern classes in $M_H$ are 
represented  by  similar sums over all faces of $\De^H$ that are 
transverse to the fiber.    Now the element $\p_i$ in the algebra 
$R_H: = \R[\p_1,\dots,\p_{N+2}]/I(V^H)$ represents the 
Poincar\'e dual to $F^{H}_i$, the prolongation of $F_i$ to $\De^H$.   
Hence the operator in $R_H$ that represents $c^{\Vertt}_s$ 
is  $\sum_{|I|=s,I\subset I_0} \p_I$, where 
$I_0: =  \{1,\dots,N\}$.

Therefore we must show that  $H$ is fully mass linear if and only if 
$$
\Psi\Bigl(\sum_{|I|=s,I\subset I_0} \p_I\Bigr) - \sum_{|I|=s,I\subset I_0} \p_I \in I(V^H),\;\; s=1,\dots,n.
$$
For simplicity, in the sums below we assume without explicit mention that
 $I\subset I_0$.  Then  we have 
\begin{eqnarray*}
\Psi\Bigl(\sum_{|I|=s} \p_I\Bigr)V^H &=& \sum_{|I|=s} \p'_I V^H\\
&=& \sum_{|I|=s}\; \prod_{i\in I} \bigl(\p_i + \ga_i \p_{N+1}\bigr)V^H
\\
&=& \sum_{|I|=s}\p_I V^H + \sum_{|I|=s, i\in I} 2\ga_i \p_{I\less i}V.
\end{eqnarray*}
We saw in the proof of Lemma \ref{le:FML} that the vanishing of the second sum above is equivalent to the identity $H(B_{n-s}) = H(B_n)$.
 Therefore (ii) holds.  Moreover (iii) holds 
 by Proposition \ref{prop:FML}. 
\end{proof}

The next result concerns the homomorphisms $I_\be$
of equation \eqref{eq:be0}.

\begin{cor}\labell{cor:cohom} Lemma \ref{le:tot} holds.
\end{cor}
\begin{proof}  Suppose that $H$ is fully mass linear.   
Since $c^M_\be$ is a product of the classes $c_i^\be$ and
$c_\be^\Vertt$ is the corresponding product of  the $c_i^\Vertt$,
 the isomorphism $\Psi$ above takes $c^M_\be$ to the class 
$c_\be^\Vertt$ for all $\be$.  
Moreover,
because the coupling class $u$ is the unique extension of $[\om]$ such that 
$u^{n+1}=0$, $\Psi$ takes the coupling class of the product to that for $M^H$.
Hence we can evaluate the integral $I_\be(H)$ of \eqref{eq:be0} on the product, where it vanishes.
\end{proof}
%\begin{rmk}\rm  There is an analogous interpretation of the totally mass linear  condition  discussed in Remark \ref{rmk:TML}.   Note that this might give new information because the vertical Chern classes may not be closed under multiplication; for example  $(c_1^{\Vertt})^2$ need not equal 
%$(c_1^2)^{\Vertt}.$
%\end{rmk}

\end{document}